\documentclass{elsarticle}

\usepackage{amsmath,amsfonts,amssymb,amsthm}
\usepackage{enumitem}
\usepackage{geometry}[margin=2cm]
\newtheorem{theorem}{Theorem}[section]
\newdefinition{definition}[theorem]{Definition}
\newtheorem{proposition}[theorem]{Proposition}
\newtheorem{corollary}[theorem]{Corollary}
\newtheorem{lemma}[theorem]{Lemma}
\newtheorem{remark}[theorem]{Remark}
\allowdisplaybreaks

\begin{document}

\title{\textbf{Stochastic integration in Riemannian manifolds from a functional-analytic point of view}}
\author{Alexandru Mustăţea}
\ead{Alexandru.Mustatea@imar.ro}
\address{"Simion Stoilow" Institute of Mathematics of the Romanian Academy, \\ P.O. Box 1-764, Bucharest, RO-014700, Romania}
\date{}

\begin{abstract}
This article presents a construction of the concept of stochastic integration in Riemannian manifolds from a purely functional-analytic point of view. We show that there are infinitely many such integrals, and that any two of them are related by a simple formula. We also find that the Stratonovich and Itô integrals known to probability theorists are two instances of the general concept constructed herein.
\end{abstract}

\begin{keyword}
stochastic integral \sep Itô integral \sep Stratonovich integral \sep Wiener measure \sep Riemannian manifold

\MSC[2020] 60H05 \sep 58J90 \sep 58J65 \sep 46E30 \sep 35K08 \sep 28C20
\end{keyword}

\maketitle

\section{Motivation and context}

The concept of stochastic integral is familiar to most probability theorists, manifesting itself in the guise of its two avatars: the Stratonovich integral and the Itô integral; it is always presented within the conceptual framework of probability theory. The aim of this work is to reconstruct the very same concept solely upon functional-analytic and Riemannian foundations. Not only shall we achieve this goal, but we shall even be able to exhibit an infinite family of such integrals, all of them particular instances of a single underlying general concept; among these we shall also find the two historically important integrals mentioned above.

In the following, $M$ will be a separable connected Riemannian manifold and $x_0 \in M$ some fixed arbitrary point. If $t>0$, we shall repeatedly make use of the space $\mathcal C _t = \{ c:[0,t] \to M \mid c \text{ is continuous, with } c(0) = x_0 \}$, that we shall endow with the natural Wiener measure $w_t$. The form to integrate along curves will be $\alpha \in \Omega^1 (M)$, a real smooth $1$-form. If $c : [0,t] \to M$ is a smooth curve, we know how to give a meaning to the line integral $\int _c \alpha$.

In order to connect this article with the stochastic literature, let us briefly recall some elements of stochastic integration in $\mathbb R^n$ without any claim of rigour. If $c$ is a smooth enough curve, the Riemann sums
\[ \sum _{j=0} ^{2^k-1} \alpha \left( c \left( \frac {jt} {2^k} \right) \right) \left[ c \left( \frac {(j+1)t} {2^k} \right) - c \left( \frac {jt} {2^k} \right) \right] \]
converge to the line integral $\int _c \alpha$. It is worth asking ourselves: if $c$ is merely continuous (or, even less, only an element of $\prod _{s \in [0,t]} M$), do these sums still converge to something meaningful and useful? The answer is known to be in the affirmative, but in a slightly weaker sense, it no longer being true for every curve: it turns out that the limit still exists, but only in measure (with respect to the Wiener measure); it is called the \textit{Itô integral} of $\alpha$. Furthermore, if we symmetrize the above Riemann sums, meaning that we should now consider the sums
\[ \sum _{j=0} ^{2^k-1} \frac 1 2 \left[ \alpha \left( c \left( \frac {jt} {2^k} \right) \right) + \alpha \left( c \left( \frac {(j+1)t} {2^k} \right) \right) \right] \left[ c \left( \frac {(j+1)t} {2^k} \right) - c \left( \frac {jt} {2^k} \right) \right] \ , \]
these, too, will converge in measure, but this time to a different limit, called the \textit{Stratonovich integral} of $\alpha$.

The starting point of our development is the useful remark that if in the formula
\[ \sum _{j=0} ^{2^k-1} \int _{[0,1]} \alpha_{(1-\tau) \, c(\frac {jt} {2^k}) + \tau \, c(\frac {(j+1)t} {2^k})} \left[ c \left( \frac {(j+1)t} {2^k} \right) - c \left( \frac {jt} {2^k} \right) \right] \mathrm d P (\tau) \]
we take the Borel probability $P$ on $[0,1]$ to be either $P = \delta_0$ (the Dirac probability concentrated at $0$) or $P = \frac 1 2 (\delta_0 + \delta_1)$, we obtain precisely the sums seen above that converge to either the Itô or, respectively, the Stratonovich integral. We conclude that these two stochastic integrals and their approximating sums seem to be particular cases of a general, single concept, that we shall indeed construct below. The generalization of this formula from $\mathbb R^n$ to $M$ is quite straightforward: the line segment $\tau \mapsto \tau \, c(\frac {jt} {2^k}) + (1-\tau) \, c(\frac {(j+1)t} {2^k})$ will get replaced by the unique minimizing geodesic between $c(\frac {jt} {2^k})$ and $c(\frac {(j+1)t} {2^k})$ (whenever it exists, of course), and the vector $c(\frac {(j+1)t} {2^k}) - c(\frac {jt} {2^k})$ will get replaced by the tangent vector to this geodesic at $\tau$.

Let us consider the trivial vector bundle $M \times \mathbb C$, endowed with the usual Hermitian structure, and with the connection $\nabla^{(\alpha)} f = \mathrm d f + \mathrm i f \alpha$, where $\mathrm i = \sqrt {-1}$ is a complex square root of $-1$. It is easy to see that $\nabla ^{(\alpha)}$ is Hermitian, and that the operator $-\Delta^{(\alpha)} = (\nabla ^{(\alpha)})^* \nabla ^{(\alpha)} : C_0 ^\infty(M) \to C_0 ^\infty(M)$ is symmetric and positive-definite. The usual Friedrichs construction will then give us a self-adjoint and positive-definite extension $H_\alpha$ that will be densely defined in $L^2(M)$. Using the results obtained by Batu Güneysu in chapter XI of his monograph \cite{Guneysu17}, the semigroup $(\mathrm e ^{-t H_\alpha})_{t \ge 0}$ will admit an integral kernel $h_\alpha$. Using the main theorem in \cite{Mizohata57} on $(0,\infty) \times M \times M$, the parabolic operator $\partial_t + H_\alpha \oplus H_\alpha$ will be hypoelliptic, whence we deduce that $h_\alpha$ is smooth. The diamagnetic inequality (proposition XI.5 in \cite{Guneysu17}), then tells us that $|h_\alpha (t,x,y)| \le h(t,x,y)$ for every $t>0$ and $x,y \in M$, where $h$ is the heat kernel on $M$.

For every $k \in \mathbb N$ we shall consider the natural projection $\pi_k : \mathcal C_t \to M ^{2^k}$ given by $\pi_k (c) = \big( c(\frac t {2^k}), \dots, c(\frac {2^k t} {2^k}) \big)$. Regardless of whether we endow $\mathcal C_t$ with the topology of uniform convergence of curves, or with the one of pointwise convergence of curves, $\pi_k$ will be continuous.

The continuous functions on some topological space will be denoted by $C(X)$, and the continuous bounded functions by $C_b (X)$. The compactly-supported smooth functions on $M$ will be $C_0 ^\infty (M)$. The complex spaces $L^p (X)$ will have the usual meaning for $p \in [1, \infty]$ whenever $X$ is endowed with a measure. The space $L^0 (X)$ is the space of complex-valued measurable functions identified under equality almost everywhere; the natural topology upon it is the one of convergence in measure.

In order to ease the reader's navigation through the text that follows, now is the right time to sketch the result that we are looking for, and the strategy that we shall use to obtain it. We shall begin by constructing a very special function $\rho_{\alpha, t} \in L^\infty (\mathcal C_t)$, following which we shall show that the map $\mathbb R \ni s \mapsto \rho_{s \alpha, t} \in \mathcal B (L^2 (\mathcal C_t))$ (the space of bounded operators in $L^2 (\mathcal C_t)$) is a strongly continuous $1$-parameter unitary group which, by Stone's theorem, will have a self-adjoint generator $\operatorname{Strat}_t ( \alpha)$ (which will be later seen to be precisely the Stratonovich stochastic integral, this also justifying its notation). The difficulty in proving this assertion comes from the fact that $\rho_{\alpha, t}$ will be obtained through an abstract procedure which will obscure the group structure and its unitarity. In order to obtain these very concrete properties, we shall construct a sequence of functions that will trivially exhibit them, and which converges to $\rho_{\alpha, t}$; this convergence will transfer these properties to $\rho_{\alpha, t}$.

More precisely, we shall construct a sequence of real measurable functions $S_{P,t,k} (\alpha)$, linear in $\alpha \in \Omega^1 (M)$, such that $\mathrm e ^{\mathrm i S_{P,t,k} (\alpha)} \to \rho_{\alpha, t}$ in $L^2 (\mathcal C_t)$. Although simple, this idea is complicated by technical details that we shall point out when we encounter them, and that force us to approach the problem indirectly: instead of proving the desired convergence directly on $L^2 (\mathcal C_t)$ (which seems extremely difficult), we shall first prove it in the space $L^2 (\mathcal C_t (\overline U))$ associated to an arbitrary relatively compact open subset $U$ with smooth boundary, following which we shall consider an exhaustion of $M$ with such subsets, which will allow us to prove the convergence in  $L^2 (\mathcal C_t)$.

\section{A generalized Wiener measure} \label{rho}

Let $U \subseteq M$ be a connected relatively compact open subset, with (possibly empty) smooth boundary, such that $x_0 \in U$ (if $M$ is compact we shall take $U = M$). We shall endow the space
\[ \mathcal C_t (\overline U) = \{ c:[0,t] \to \overline U \mid c \text{ is continuous, with } c(0) = x_0 \} \]
with the corresponding intrinsic Wiener measure $w_{t} ^{(U)}$ (for details about the Wiener measure, the article \cite{BP11} contains all the necessary constructions and explanations; note that the constructions therein are not probabilistic, but functional-analytic, therefore our project of a purely functional-analytic construction of stochastic integration is not compromised). This is a metric space when endowed with the distance $D(c_0, c_1) = \max_{s \in [0,t]} d \big( c_0(s), c_1(s) \big)$; it is separable (and therefore second-countable) by \cite{Michael61}. In particular, we may use Luzin's theorem on it.

Let
\[ \operatorname{Cyl} (\mathcal C_t (\overline U)) = \{ f \in C_b (\mathcal C_t (\overline U)) \mid \exists k \in \mathbb N \text{ and } f_k \in C(\overline U ^{2^k}) \text{ such that } f = f_k \circ \pi_k \} \]
be the algebra of continuous cylindrical functions on $\mathcal C_t (\overline U)$. Clearly, $\operatorname{Cyl} (\mathcal C_t (\overline U)) \subset L^1 (\mathcal C_t (\overline U))$.

\begin{theorem}
The algebra $\operatorname{Cyl} (\mathcal C_t (\overline U))$ is dense in $L^p (\mathcal C_t (\overline U), w_t ^{(U)})$ for every $p \in [1, \infty)$.
\end{theorem}
\begin{proof}
The methods chosen for the proof will need the order relationship on $\mathbb R$, therefore we shall first assume that all the function spaces involved are real; the case of complex functions will then follow trivially from the real one.

The proof strategy is the following: first, we shall approximate the functions in $L^p$ with bounded functions in $L^p$; next, we shall consider a compact subset with sufficiently small complementary subset, on which we shall approximate the bounded functions in $L^p$ with continuous cylindrical functions using the Stone-Weierstrass theorem; finally, we shall show that we can control the behaviour of these approximating cylindrical functions on the complementary subset of the chosen compact subset.

That the space $L^p _b (\mathcal C_t (\overline U))$ (the essentially bounded functions in $L^p (\mathcal C_t (\overline U))$) is dense in $L^p (\mathcal C_t (\overline U))$ is obvious since the successive truncations of any function converge to it.

If now $0 \ne f \in L^p _b (\mathcal C_t (\overline U))$, let $B = \operatorname{ess \, sup} |f|$ and let $\varepsilon \in (0, \min(B,1))$. Choose $f_0 : \mathcal C_t (\overline U) \to \mathbb R$ a measurable representative of $f$ with $\sup |f_0| = B$. Using Luzin's theorem, there exists a co-null subset $\emptyset \ne \mathcal S \subseteq \mathcal C_t (\overline U)$ such that $f_0 | _{\mathcal S}$ is continuous in the topology induced on $\mathcal S$. The restriction of $w_t ^{(U)}$ to $\mathcal S$ is a non-trivial regular Borel measure; let then $\mathcal K \subseteq \mathcal S$ be a compact subset (in the induced topology) such that $(w_t ^{(U)} | _{\mathcal S}) (\mathcal S \setminus \mathcal K) < \frac \varepsilon {2 (2^p + 1) B^p}$. Let us show that the subalgebra $\operatorname{Cyl} (\mathcal C_t (\overline U)) | _{\mathcal K}$ consisting of the functions in $\operatorname{Cyl} (\mathcal C_t (\overline U))$ restricted to $\mathcal K$ satisfies the hypotheses of the real version of the Stone-Weierstrass theorem on compact spaces.

First, it is obvious that $1 = 1 \circ \pi_0 \in \operatorname{Cyl} (\mathcal C_t (\overline U)) | _{\mathcal K}$. It remains to show that $\operatorname{Cyl} (\mathcal C_t (\overline U))$ separates the points of $\mathcal C_t (\overline U)$, whence in particular it will result that $\operatorname{Cyl} (\mathcal C_t (\overline U)) | _{\mathcal K}$ separates the points of $\mathcal K$. Let then $c_0, c_1 \in \mathcal C_t (\overline U)$ with $c_0 \ne c_1$. Since these curves are continuous and the "dyadic" numbers $\{ \frac {jt} {2^k} \in (0,1] \mid k \in \mathbb N, \, 1 \le j \le 2^k \}$ are dense in $(0,t]$, there exist $k \in \mathbb N$ and $j \in \{ 1, \dots, 2^k \}$ such that $c_0 (\frac {jt} {2^k}) \ne c_1 (\frac {jt} {2^k})$. If $\varphi_j : \overline U \to [0,1]$ is a continuous function that separates the points $c_0 (\frac {jt} {2^k})$ and $c_1 (\frac {jt} {2^k})$, then the function $(1 \otimes \dots \otimes 1 \otimes \varphi_j \otimes 1 \dots \otimes 1) \circ \pi_k \in \operatorname{Cyl} (\mathcal C_t (\overline U))$ clearly separates $c_0$ and $c_1$.

We may now apply the Stone-Weierstrass theorem on $\mathcal K$, whence there exists $g' \in \operatorname{Cyl} (\mathcal C_t (\overline U))$ such that $\sup | f_0 | _{\mathcal K} - g' | _{\mathcal K} | < \frac \varepsilon 2$, whence it will follow that $\| f | _{\mathcal K} - g' | _{\mathcal K} \| _{L^p (\mathcal K, w_t ^{(U)})} < \frac \varepsilon 2$.

So far, we have obtained a cylindrical function $g'$ that approximates $f$ on $\mathcal K$. It remains to see what to do with $g'$ on $\mathcal S \setminus \mathcal K$. To this end, let $\varphi : \mathbb R \to \mathbb R$ be given by
\[ \varphi (r) = \begin{cases}
-2B, & r < -2B , \\
r, & r \in [-2B, 2B], \\
2B, & r > 2B \end{cases} \ . \]
Clearly, $\varphi$ is continuous and bounded, so the function $g = \varphi \circ g'$ (it is here where we use that $g'$ is real) belongs to $\operatorname{Cyl} (\mathcal C_t (\overline U))$. Notice that
\[ \sup | g | _{\mathcal K} | \le \sup | g | _{\mathcal K} - f_0 | _{\mathcal K} | + \sup | f_0 | _{\mathcal K} | \le \frac \varepsilon 2 + B < 2B \ , \]
so $g | _{\mathcal K} = g' | _{\mathcal K}$. Putting all these ingredients together, we obtain that
\begin{align*}
\| f-g \| _{L^p (\mathcal C_t (\overline U), w_t ^{(U)})} ^p & = \| f-g \| _{L^p (\mathcal K, w_t ^{(U)})} ^p + \| f-g \| _{L^p (\mathcal S \setminus \mathcal K, w_t ^{(U)})} ^p + \| f-g \| _{L^p (\mathcal C_t (\overline U) \setminus \mathcal S, w_t ^{(U)})} ^p = \\
& = \| f-g' \| _{L^p (\mathcal K, w_t ^{(U)})} ^p + \int _{\mathcal S \setminus \mathcal K} |f-g| ^p \, \mathrm d w_t ^{(U)} + 0 < \\
& < \left( \frac \varepsilon 2 \right) ^p + (2^p + 1) B^p \, w_t ^{(U)} (\mathcal C_t (\overline{U}) \setminus \mathcal K) < \frac \varepsilon 2 + \frac \varepsilon 2 = \varepsilon \ ,
\end{align*}
which proves that $\operatorname{Cyl} (\mathcal C_t (\overline U))$ is dense in $L^p (\mathcal C_t (\overline U))$.
\end{proof}

Let us define the (obviously linear) functional $W_{\alpha, t} ^{(U)}: \operatorname{Cyl} (\mathcal C_t(\overline U)) \to \mathbb C$ by
\[ W_{\alpha, t} ^{(U)} (f_k \circ \pi_k) = \int_U \mathrm d x_1 \, h_\alpha ^{(U)} \left( \frac t {2^k}, x_0, x_1 \right) \dots \int_U \mathrm d x_{2^k} \, h_\alpha ^{(U)} \left( \frac t {2^k}, x_{2^k - 1}, x_{2^k} \right) f_k (x_1, \dots, x_{2^k}) \]
for every $f_k \circ \pi_k \in \operatorname{Cyl} (\mathcal C_t (\overline U))$, where $h _\alpha ^{(U)}$ is the integral kernel on $U$ associated to the connection $\mathrm d + \mathrm i \alpha$ in the trivial bundle $U \times \mathbb C$, constructed as explained above (again, for details see chapter XI of \cite{Guneysu17}). The next theorem will produce a measure density on $\mathcal C_t (\overline U))$ that will depend on the form $\alpha$ and that will be the main object of study in the first half of this article. Its product with the Wiener measure may be thought of as a generalized, or perturbed, Wiener measure; when $\alpha=0$ it coincides with the usual Wiener measure.

\begin{theorem}
There exists a unique $\rho_{\alpha, t} ^{(U)} \in L^\infty (\mathcal C_t (\overline U))$ with $\| \rho_{\alpha, t} ^{(U)} \| _{L^\infty (\mathcal C_t (\overline U))} \le 1$ such that $W_{\alpha, t} ^{(U)} (f) = \int_{\mathcal C_t (\overline U)} f \rho_{\alpha, t} ^{(U)} \, \mathrm d w_t ^{(U)}$ for every $f \in \operatorname{Cyl} (\mathcal C_t (\overline U))$.
\end{theorem}
\begin{proof}
Taking the absolute value in the definition of $W_{\alpha, t} ^{(U)}$ and using the diamagnetic inequality, we have that
\begin{align*}
| W_{\alpha, t} ^{(U)} & (f_k \circ \pi_k) | \le \\
& \le \int_U \mathrm d x_1 \, \left| h_\alpha ^{(U)} \left( \frac t {2^k}, x_0, x_1 \right) \right| \dots \int_U \mathrm d x_{2^k} \, \left| h_\alpha ^{(U)} \left( \frac t {2^k}, x_{2^k - 1}, x_{2^k} \right) \right| |f_k (x_1, \dots, x_{2^k})| \le \\
& \le \int_U \mathrm d x_1 \, h ^{(U)} \left( \frac t {2^k}, x_0, x_1 \right) \dots \int_U \mathrm d x_{2^k} \, h ^{(U)} \left( \frac t {2^k}, x_{2^k - 1}, x_{2^k} \right) |f_k (x_1, \dots, x_{2^k})| = \\
& = \| f_k \circ \pi_k \| _{L^1 (\mathcal C_t (\overline U))} \ ,
\end{align*}
so $W_{\alpha, t} ^{(U)}$ is continuous in the norm $\| \cdot \| _{L^1(\mathcal C_t (\overline U))}$ on $\operatorname{Cyl} (\mathcal C_t (\overline U))$; since the latter is dense in $L^1(\mathcal C_t (\overline U))$, it follows that we may extend $W_{\alpha, t} ^{(U)}$ to a continuous linear functional on $L^1(\mathcal C_t (\overline U))$, hence there exists $\rho_{\alpha, t} ^{(U)} \in L^\infty (\mathcal C_t (\overline U))$ such that $W_{\alpha, t} ^{(U)} (f) = \int_{\mathcal C_t (\overline U)} f \rho_{\alpha, t} ^{(U)} \, \mathrm d w_t ^{(U)}$ for every $f \in L^1 (\mathcal C_t (\overline U))$. Furthermore, $|W_{\alpha, t} ^{(U)} (f)| \le \| f \| _{L^1(\mathcal C_t (\overline U))}$ for every $f \in L^1 (\mathcal C_t (\overline U))$, so $\| \rho_{\alpha, t} ^{(U)} \| _{L^\infty (\mathcal C_t (\overline U))} \le 1$.
\end{proof}

\section{A sequence of approximations for $\rho_{\alpha, t} ^{(U)}$}

So far, $\rho_{\alpha, t} ^{(U)}$ has been constructed by a very abstract argument, therefore its various concrete properties are difficult to study. As a consequence, in what follows we shall construct a sequence of concrete approximations of this function, which will enjoy two essential properties: a group property, and the fact of being of absolute value $1$. We shall then show that this sequence converges to $\rho_{\alpha, t} ^{(U)}$ in $L^2 (\mathcal C_t)$, so that these two properties will be transferred to $\rho_{\alpha, t} ^{(U)}$, too. In order to complete this program, we shall now introduce several more ingredients.

Let $P$ be a Borel regular probability on $[0,1]$; we shall see later on that the role of $P$ will be to classify the various stochastic integrals that we shall obtain. Let $M_1(P)$ be the first order moment of $P$, that is
\[M_1(P) = \int _{[0,1]} \tau \, \mathrm d P (\tau) \ . \]

Whenever the points $x, y \in M$ may be joined by a unique minimizing geodesic, we shall denote it by $\gamma_{x,y} : [0,1] \to M$, where we understand that $\gamma(0) = x$ and $\gamma(1) = y$. Let us now define $I_P (\alpha) : M \times M \to \mathbb R$ by:
\begin{itemize}
\item $I_P (\alpha) (x,y) = \int _{[0,1]} \alpha _{\gamma_{x,y} (\tau)} (\dot \gamma_{x,y} (\tau)) \, \mathrm d P (\tau)$, if there exists a unique minimizing geodesic $\gamma_{x,y}$ as above between $x$ and $y$;
\item $I_P (\alpha) (x,y) = 0$, otherwise.
\end{itemize}

\begin{proposition}
$I_P (\alpha)$ is smooth on $\{ (x,y) \in M \times M \mid x \in M, \, y \in \operatorname{Dom} (\exp_x) \}$, where $\operatorname{Dom} (\exp_x)$ is the domain of definition of the Riemannian exponential $\exp_x$ at $x$.
\end{proposition}
\begin{proof}
The proof being elementary, we shall only sketch it. On $\{ (x,y) \in M \times M \mid x \in M, \, y \in \operatorname{Dom} (\exp_x) \}$ $I_P (\alpha)$ may be written explicitly as
\[ I_P (\alpha) (x,y) = \int _{[0,1]} \alpha _{\exp_x (\tau \exp_x ^{-1} (y))} \, [(\mathrm d \exp_x) _{\tau \exp_x ^{-1} (y)} (\exp_x ^{-1} (y))] \, \mathrm d P (\tau) \ ,\]
which is seen to be smooth from the smoothness of the Riemannian exponential in both arguments, followed by an application of the dominated convergence theorem.
\end{proof}

The next ingredient to introduce will be a smooth cut-off function $\chi$, the role of which being to keep us away from the points where $I_P (\alpha)$ stops being smooth. To this end, let $\kappa : [0, \infty) \to [0,1]$ be a smooth function such that $\kappa | _{[0, \frac 1 3]} = 1$ and $\kappa | _{[\frac 1 2, \infty)} = 0$. Let $\operatorname{injrad}_U : U \to (0, \infty)$ be the injectivity radius function on $U$; we emphasize that this is not the restriction of $\operatorname{injrad}_M$ to $U$, but rather it is computed intrinsically, using the restriction to $U$ of the Riemannian structure (for basic details about the injectivity radius, see p.118 of \cite{Chavel06}). Being continuous and strictly positive, we may find a smooth function $r : U \to (0, \infty)$ such that $r(x) < \operatorname{injrad}_U (x)$. In particular, $r(x) \le d_U (x, \partial U)$ (the distance up to the boundary of $U$, computed using the intrinsic distance $d_U$ of $U$, not using the distance on $M$ restricted to $U$). We may now finally define the desired cut-off function $\chi : U \times U \to [0,1]$ by $\chi (x,y) = \kappa \left( \frac {d_U (x,y)^2} {r(x) ^2} \right)$. Notice that $\chi$ is smooth (the square is necessary in order to guarantee the smoothness close to the points with $y=x$).

All the ingredients introduced so far in this subsection were necessary in order for us to be able to construct the operator $R_{\alpha, t} ^{(U)}$ by $R_{\alpha, 0} ^{(U)} f = f$ and
\[ (R_{\alpha, t} ^{(U)} f) (x) = \int _U h_\alpha ^{(U)} (t,x,y) \, \chi(x,y) \, \mathrm e ^{-\mathrm i I_P (\alpha) (x,y) - \mathrm i t (\mathrm d^* \alpha) (x) \int _{[0,1]} (2 \tau - 1) \, \mathrm d P (\tau)} f(y) \, \mathrm d y \]
for every $t>0$ and $f \in C_b(U)$, where $\mathrm d^*$ is the Hodge codifferential, defined as the formal adjoint of the differential operator $\mathrm d$, that is
\[ \int _M (\mathrm d^* \alpha) (x) \, \varphi (x) \, \mathrm d x = \int _M \langle \alpha_x, \mathrm d_x \varphi \rangle _{T^* _x M} \, \mathrm d x \]
for every $\varphi \in C_0 ^\infty (M, \mathbb R)$ (the notation $\langle \cdot, - \rangle _{T^* _x M}$ denoting the scalar product on $T^* _x M$ induced by the Riemannian structure). For a more "Riemannian" understanding, if $\alpha^\sharp$ is the vector field dual to $\alpha$ by raising the latter's indices, then $\mathrm d^* \alpha = -\operatorname{div} \alpha^\sharp$. It is worth noting that the operator $R_{\alpha, t} ^{(U)}$ was constructed as a "geometrically hybrid" object, and it was purposefully so: on the one hand, the factors $h_\alpha ^{(U)}$ and $\chi$ stem from the intrinsic Riemannian geometry of $\overline U$; on the other hand, the factor containing $I_P (\alpha)$ is extrinsic, meaning that these functions are the restrictions to $\overline U$ of functions defined on the whole manifold $M$. We could have used intrinsic versions of $I_P (\alpha)$, but this would have made some theorems that we shall encounter next much more difficult to prove. Since the Hodge codifferential is a local operator, in its case the distinction between intrinsic and extrinsic does not really matter.

Let us notice that the integrand in the formula of $R_{\alpha, t} ^{(U)}$ is smooth: even though $I_P (\alpha)$ is, in principle, discontinuous at the points $(x,y)$ with $d(x,y) = \operatorname{injrad}_M (x)$, these points are contained in the complementary subset of the support of the cut-off function $\chi$ (because $\operatorname{injrad}_U (x) \le \operatorname{injrad}_M (x)$), and $\chi$ is smooth.

If $h^{(U)}$ is the intrinsic heat kernel of $\overline U$, the operators defined by
\[ C(\overline U) \ni f \mapsto \int _U h^{(U)} (s, \cdot ,y) \, f(y) \, \mathrm d y \in C(\overline U) \]
together with the identity operator form a strongly continuous one-parameter semigroup in $C(\overline U)$. This will have a generator (closed operator) that we shall denote by $L^{(U)}$, densely defined, with the domain given by (see \cite{Davies80}, chap. 1)
\[ \operatorname{Dom} (L^{(U)}) = \left\{ f \in C(\overline U) \mid \lim _{s \to 0} \frac 1 s \left( \int _U h^{(U)} (s, \cdot ,y) \, f(y) \, \mathrm d y - f \right) \in C(\overline U) \right\} \ . \]
We shall denote this semigroup by $(\mathrm e ^{-s L^{(U)}}) _{s \ge 0}$. An essential domain for $L^{(U)}$ is
\[ \mathcal E = \bigcup _{s > 0} \mathrm e ^{-s L^{(U)}} (C(\overline U)) \ . \]
Since the heat semigroup is smoothing (again, one may use \cite{Mizohata57}, or one's favourite Sobolev spaces techniques, to see this), the functions in $\mathcal E$ will be smooth. Since $h^{(U)}$ vanishes on the boundary $\partial U$, the functions in $\mathcal E$ will also vanish on $\partial U$.

With exactly the same arguments, but using now the integral kernel $h_\alpha ^{(U)}$ instead of $h ^{(U)}$, we shall obtain another semigroup acting on $C(\overline U)$, the generator of which will be denoted by $L_\alpha ^{(U)}$.

\begin{lemma} \label{the compactly-supported functions are in the domain}
$\operatorname{Dom} (L_\alpha ^{(U)})$ contains the space $C _0 ^\infty (U)$ of the compactly-supported smooth functions.
\end{lemma}
\begin{proof}
If $u \in C _0 ^\infty (U)$ it is clear that $L_\alpha ^{(U)} u \in C(\overline U)$. We shall show that
\[ \lim _{s \to 0} \frac 1 s (\mathrm e ^{-s L_\alpha ^{(U)}} u - u) = -L_\alpha ^{(U)} u \]
in the norm topology of $C (\overline U)$.

To begin with, let us show that $[0, \infty) \ni s \mapsto (\mathrm e ^{-s L_\alpha ^{(U)}} u) (x) \in \mathbb C$ is smooth for all $x \in \overline U$. Since $h_\alpha ^{(U)}$ is smooth, the function $(s,x) \mapsto h_\alpha ^{(U)} (s,x,y) \, u(y)$ will be smooth for all $y \in \overline U$. We shall denote by $(L_\alpha ^{(U)})_y$ the operator $L_\alpha ^{(U)}$ acting with respect to the argument $y$; since $\partial_s h(s,x,y) = -(L_\alpha ^{(U)})_y \, h(s,x,y)$ and $\overline U$ is compact, we may use the dominated convergence theorem to differentiate with respect to $s \in (0, \infty)$ under the integral sign and obtain that
\[ \frac {\mathrm d} {\mathrm d s} (\mathrm e ^{-s L_\alpha ^{(U)}} u) (x) = \int _U -(L_\alpha ^{(U)})_y \, h_\alpha ^{(U)} (s,x,y) \, u(y) \, \mathrm d y = \int _U h_\alpha ^{(U)} (s,x,y) \, (-L_\alpha ^{(U)} u) (y) \, \mathrm d y \ . \]
This argument may be iterated indefinitely, so $(0, \infty) \ni s \mapsto (\mathrm e ^{-s L_\alpha ^{(U)}} u) (x) \in \mathbb C$ is smooth for all $x \in \overline U$. Passing to the limit when $s \to 0$ also gets us the smoothness at $0$.

It is easy to see that
\begin{align*}
\lim _{s \to 0} \frac 1 s (\mathrm e ^{-s L_\alpha ^{(U)}} u - u) (x) & = \lim _{s \to 0} \partial_s \int _U h_\alpha ^{(U)} (s,x,y) \, u(y) \, \mathrm d y = \\
& = \lim _{s \to 0} \int _U h_\alpha ^{(U)} (s,x,y) \, (-L_\alpha ^{(U)} u) (y) \, \mathrm d y = (-L_\alpha ^{(U)} u) (x)
\end{align*}
for all $x \in \overline U$.

Consider now the function $F_u : [0, \infty) \to C(\overline U)$ given by $F_u (s) = \mathrm e ^{-s L_\alpha ^{(U)}} u - u + s L_\alpha ^{(U)} u$. We have that $F_u (0) (x) = 0$ and $F_u ' (0) (x) = 0$ for all $x \in \overline U$, whence it follows that
\begin{align*}
\| F_u (s) \| _{C(\overline U)} & = \sup_{x \in \overline U} | F_u (s) (x) | = \sup_{x \in \overline U} \left| \int _0 ^s (s - \sigma) F_u '' (\sigma) (x) \, \mathrm d \sigma \right| \le \\
& \le \frac {s^2} 2 \sup_{x \in \overline U} \sup_{\sigma \in [0,s]} | F_u '' (\sigma) (x) | \le \frac {s^2} 2 \sup_{x \in \overline U} \sup_{\sigma \in [0,s]} \int _{\overline U} |h_\alpha ^{(U)} (\sigma ,x,y)| \, |(L_\alpha ^{(U)})^2 u| (y) \, \mathrm d y \le \\
& \le \frac {s^2} 2 \| (L_\alpha ^{(U)})^2 u \| _{C(\overline U)} \ ,
\end{align*}
where at the jump between the second and the third line we have used the diamagnetic inequality and the sub-Markovianity of $h ^{(U)}$. This shows that $\lim _{s \to 0} \| F_u (s) \| _{C(\overline U)} = 0$, which means that $\lim _{s \to 0} \frac 1 s (\mathrm e ^{-s L_\alpha ^{(U)}} u - u) = -L_\alpha ^{(U)} u$
in the norm topology of $C(\overline U)$, as desired, whence $u \in \operatorname{Dom} L_\alpha ^{(U)}$ as claimed.
\end{proof}

The crucial tool to be used in the following will be Chernoff's theorem (lemma 3.28 in \cite{Davies80}). For the reader's convenience, we shall give its statement here.

\begin{theorem}[Chernoff]
Assume that $(R_t) _{t \ge 0}$ is a family of contractions in a Banach space $X$, with $R_0 = \operatorname{Id}_X$. Let $\mathcal E \subseteq X$ be an essential domain for the generator $L$ of a strongly continuous one-parameter semigroup $(\mathrm e ^{-t L}) _{t \ge 0}$ on $X$. If $\lim _{t \to 0} \frac 1 t (R_t f - f) = -L f$ for every $f \in \mathcal E$, then $\mathrm e ^{-t L} = \lim _{k \to \infty} \big( R_{\frac t k} \big) ^k$ strongly for every $t \ge 0$. Furthermore, the convergence is uniform with respect to $t$ in bounded subsets of $[0, \infty)$.
\end{theorem}

With all these preparations, we are ready now for the main technical result of this work, from which all the conceptual developments announced in the introduction will unravel (the author apologizes in advance for the lengthy calculations involved).

\begin{theorem} \label{the technical theorem}
For every $t \ge 0$ and for every $f \in C(\overline U)$, $\mathrm e^{-tL^{(U)}} f = \lim _{k \to \infty} \big( R_{\alpha, \frac t k} ^{(U)} \big) ^k f$ in the topology of uniform convergence on $U$, uniformly with respect to $t$ in bounded subsets of $[0, \infty)$.
\end{theorem}
\begin{proof}
The proof reduces to the verification of the assumptions in Chernoff's theorem, whence the conclusion will be immediate.

To begin with, let us notice that $R_{\alpha, t} ^{(U)}$ is a contraction: indeed,
\begin{align*}
\| R_{\alpha, t} ^{(U)} f \| & _{C(\overline U)} = \\
& = \sup _{x \in U} \left| \int _U h_\alpha ^{(U)} (t,x,y) \, \chi(x,y) \, \mathrm e ^{-\mathrm i I_P (\alpha) (x,y) - \mathrm i t (\mathrm d^* \alpha) (x) \int _{[0,1]} (2 \tau - 1) \, \mathrm d P (\tau)} \, f(y) \, \mathrm d y \right| \mathrm d x \le \\
& \le \sup_{x \in U} \int _U h^{(U)} (t,x,y) \, |f(y)| \, \mathrm d y \le \| f \|_{C(\overline U)} \ ,
\end{align*}
where we have used the diamagnetic inequality $|h_\alpha ^{(U)}| \le h^{(U)}$ and the fact that $\int _U h^{(U)} (t,x,y) \, \mathrm d y \le 1$ (sub-Markovianity).

It remains to be shown that $\lim _{t \to 0} \| \frac 1 t (R_{\alpha, t} ^{(U)} f - f) + L^{(U)} f \| _{C(\overline U)} = 0$ for every $f \in \mathcal E$ (which is an essential domain for $L^{(U)}$); to this end, let us show first that $(R_{\alpha, t} ^{(U)} f) (x)$ is smooth with respect to $t$ for every $x \in U$. We notice that
\[ (R_{\alpha, t} ^{(U)} f) (x) = \langle \delta_x, \mathrm e ^{-t L_\alpha ^{(U)}} [\chi (x, \cdot) \, \mathrm e ^{-\mathrm i I_P (\alpha) (x, \cdot) - \mathrm i t (\mathrm d^* \alpha) (x) \int _{[0,1]} (2 \tau - 1) \, \mathrm d P (\tau)} \, f] \rangle \ , \]
where $\delta_x$ is the Dirac measure concentrated at $x$ and $\langle \cdot, - \rangle$ denotes the dual pairing between the space of the complex Borel regular measures on $\overline U$ and the space $C(\overline U)$. By the construction of $\chi$, the function
\[ \chi (x, \cdot) \, \mathrm e ^{-\mathrm i I_P (\alpha) (x, \cdot) - \mathrm i t (\mathrm d^* \alpha) (x) \int _{[0,1]} (2 \tau - 1) \, \mathrm d P (\tau)} \]
is smooth with compact support (as mentioned above, the possible singularities of $I_P (\alpha) (x, \cdot)$ live away from the support of $\chi(x, \cdot)$), and $f$ is smooth since it belongs to $\mathcal E$, so their product is a smooth function with compact support, therefore it belongs to the domain of every power of $L_\alpha ^{(U)}$, with the argument in lemma \ref{the compactly-supported functions are in the domain}. Under these circumstances, the map
\[ [0, \infty) \ni t \mapsto \mathrm e ^{-t L^{(U)} _\alpha} [\chi (x, \cdot) \, \mathrm e ^{-\mathrm i I_P (\alpha) (x, \cdot) - \mathrm i t (\mathrm d^* \alpha) (x) \int _{[0,1]} (2 \tau - 1) \, \mathrm d P (\tau)} \, f] \in C(\overline U) \]
is smooth, therefore it will remain so when the measure $\delta_x$ seen as an element from $C(\overline U) ^*$ is applied to it.

Considering the Taylor expansion of $(R_{\alpha, t} ^{(U)} f) (x)$ with respect to $t$ we have, for every $x \in \overline U$,
\begin{equation}
(R_{\alpha, t} ^{(U)} f) (x) = f(x) + \partial_t | _{t=0} (R_{\alpha, t} ^{(U)} f) (x) \, t + \int _0 ^t (t-s) \, \partial_s ^2 (R_{\alpha, s} ^{(U)} f) (x) \, \mathrm d s \ . \label{Taylor expansion}
\end{equation}

For the computation of the first derivative $(R_{\alpha, t} ^{(U)} f) (x)$ we have
\begin{align}
\nonumber \partial_t | _{t=0} & [(R_{\alpha, t} ^{(U)} f) (x)] = \\
\nonumber & = \lim _{t \to 0} \int _U [-L_{\alpha, y} ^{(U)} h_\alpha ^{(U)} (t,x,y)] \, \chi(x,y) \, \mathrm e ^{-\mathrm i I_P (\alpha) (x,y)} f(y) \, \mathrm d y - \\
\nonumber & - \lim _{t \to 0} \int _U h_\alpha ^{(U)} (t,x,y) \, \mathrm e ^{-\mathrm i I_P (\alpha) (x,y)} \chi(x,y) \left[ \mathrm i (\mathrm d^* \alpha) (x) \int _{[0,1]} (2 \tau - 1) \, \mathrm d P (\tau) \right] f(y) \, \mathrm d y = \\
\nonumber & = -L_\alpha ^{(U)} [\chi(x, \cdot) \, \mathrm e ^{-\mathrm i I_P (\alpha) (x,\cdot)} f] (x) - \left[ \mathrm i (\mathrm d^* \alpha) (x) \int _{[0,1]} (2 \tau - 1) \, \mathrm d P (\tau) \right] f(x) = \\
& = \Delta^{(\alpha)} [\chi(x, \cdot) \, \mathrm e ^{-\mathrm i I_P (\alpha) (x,\cdot)} f] (x) - \left[ \mathrm i (\mathrm d^* \alpha) (x) \int _{[0,1]} (2 \tau - 1) \, \mathrm d P (\tau) \right] f(x) \ , \label{Chernoff condition}
\end{align}
because $L_\alpha ^{(U)} = -\Delta^{(\alpha)}$ on functions from $C _0 ^\infty (U)$, $I_P (\alpha) (x,x) = 0$ and $\chi(x, \cdot) = 1$ in some neighbourhood of $x$ by construction.

In order to evaluate the term that contains $\Delta^{(\alpha)}$ we shall work in normal coordinates around $x$ and use the formula
\[\Delta^{(\alpha)} = -(\nabla^{(\alpha)})^* \nabla^{(\alpha)} = \sum_{i, j = 1} ^n g^{ij} \nabla^{(\alpha)} _{\partial_i} \nabla^{(\alpha)} _{\partial_j} - \sum_{i, j, k = 1} ^n g^{ij} \Gamma_{ij} ^k \nabla^{(\alpha)} _{\partial_k} \ , \]
whence, for every complex-valued smooth function $F$,
\begin{align}
\nonumber [\Delta^{(\alpha)} & F] (x) = \sum_{j=1} ^n [\nabla^{(\alpha)} _{\partial_j} \nabla^{(\alpha)} _{\partial_j} F] (x) = \sum_{j=1} ^n (\partial_j + \mathrm i \alpha_j) (\partial_j F + \mathrm i \alpha_j F) (x) = \\
\nonumber & = \sum_{j=1} ^n (\partial_j ^2 F) (x) + \mathrm i \sum_{j=1} ^n ( \partial_j \alpha_j ) (x) F (x) + 2 \mathrm i \sum_{j=1} ^n \alpha_j (x) (\partial_j F) (x) - \sum_{j=1} ^n \alpha_j (x) ^2 F(x) = \\
& = (\Delta F) (x) - \mathrm i (\mathrm d ^* \alpha) (x) F(x) + 2 \mathrm i \langle \alpha_x, \mathrm d _x F \rangle _{T^* _x M \otimes \mathbb C} - \| \alpha_x \| _{T^* _x M} ^2 F(x) \ , \label{laplacian with form}
\end{align}
where $\langle \alpha_x, \mathrm d _x F \rangle _{T^* _x M \otimes \mathbb C}$ is simply $\langle \alpha_x, \mathrm d _x u \rangle _{T^* _x M} + \mathrm i \langle \alpha_x, \mathrm d _x v \rangle _{T^* _x M} $ when $F = u + \mathrm i v$.

We shall need the explicit expressions of $[\Delta I_P (\alpha) (x, \cdot)] (x)$ and $\mathrm d _x I_P (\alpha) (x, \cdot)$ below, so we shall take this opportunity to derive them now. In the following, $\Delta_{(y)}$ and $\mathrm d _{(y)}$ will denote the Laplacian and the differential computed with respect to $y \in U$. In normal coordinates on $U$ centered at $x$ (in which, for notational simplicity, we shall identify the points around $x$ with their images in $T_x M$ under the inverse of the Riemannian exponential map $\exp_x ^{-1}$ at $x$), the geodesic $\gamma_{x,y} (\tau)$ becomes the line segment $[0,1] \ni \tau \mapsto x + \tau (y-x)$, therefore
\begin{align*}
[\Delta_{(y)} & I_P (\alpha) (x, y)] | _{y=x} = \sum _{j=1} ^n \sum _{k=1} ^n \int _{[0,1]} \partial_j ^2 [\alpha_k (x + \tau (y-x)) (y_k-x_k)] | _{y=x} \, \mathrm d P (\tau)  = \\
& = 2 \sum _{j=1} ^n \int _{[0,1]} \partial_j [\alpha_j (x + \tau (y-x))] | _{y=x} \, \mathrm d P (\tau)  = 2 \sum _{j=1} ^n \int _{[0,1]} (\partial_j \alpha_j) (x) \, \tau \, \mathrm d P (\tau) = \\
& = - 2 (\mathrm d ^* \alpha) (x) \int _{[0,1]} \tau \, \mathrm d P (\tau) \ ,
\end{align*}
and
\[ [\mathrm d _{(y)} I_P (\alpha) (x, y)] | _{y=x} = \int_{[0,1]} \mathrm d _{(y)} [\alpha_{x + \tau (y-x)} (y-x)] | _{y=x} \, \mathrm d P (\tau) = \alpha_x \ . \]

Taking $F = \mathrm e ^{-\mathrm i I_P (\alpha) (x,\cdot)} f$ we have on the one hand
\begin{align*}
\Delta [\mathrm e ^{-\mathrm i I_P (\alpha) (x,\cdot)} & f] (x) = [\Delta \mathrm e ^{-\mathrm i I_P (\alpha) (x,\cdot)}] (x) f(x) + 2 \langle \mathrm d _x \mathrm e ^{-\mathrm i I_P (\alpha) (x,\cdot)}, \mathrm d _x f \rangle _{T^* _x M \otimes \mathbb C} + (\Delta f) (x) = \\
& = - \mathrm i [\Delta I_P (\alpha) (x, \cdot)] (x) f(x) - \| \mathrm d _x I_P (\alpha) (x, \cdot) \| _{T^* _x M} ^2 f(x) - \\
& - 2 \mathrm i \langle \mathrm d _x I_P (\alpha) (x, \cdot), \mathrm d _x f \rangle _{T^* _x M \otimes \mathbb C} + (\Delta f) (x) = \\
& = 2 \mathrm i (\mathrm d^* \alpha) (x) f(x) - \| \alpha_x \| _{T^* _x M} ^2 f(x) - 2 \mathrm i \langle \alpha_x, \mathrm d _x f \rangle _{T^* _x M \otimes \mathbb C} + (\Delta f) (x)
\end{align*}
and on the other hand
\begin{align*}
\langle \alpha_x, \mathrm d _x [\mathrm e ^{-\mathrm i I_P (\alpha) (x,\cdot)} f] \rangle _{T^* _x M \otimes \mathbb C} & = \langle \alpha_x, \mathrm d _x [-\mathrm i I_P (\alpha) (x,\cdot)] \rangle _{T^* _x M \otimes \mathbb C} f(x) + \langle \alpha_x, \mathrm d _x f \rangle _{T^* _x M \otimes \mathbb C} = \\
& = - \mathrm i \| \alpha_x \| _{T^* _x M} ^2 f(x) + \langle \alpha_x, \mathrm d _x f \rangle _{T^* _x M \otimes \mathbb C} \ ,
\end{align*}
whence, returning to formula \eqref{laplacian with form},
\begin{align*}
\Delta & ^{(\alpha)} [\mathrm e ^{-\mathrm i I_P (\alpha) (x,\cdot)} f] (x) = \\
& = 2 \mathrm i (\mathrm d^* \alpha) (x) \int _{[0,1]} \tau \, \mathrm d P (\tau) \, f(x) - \| \alpha_x \| _{T^* _x M} ^2 f(x) - 2 \mathrm i \langle \alpha_x, \mathrm d _x f \rangle _{T^* _x M \otimes \mathbb C} + \\
& + (\Delta f) (x) - \mathrm i (\mathrm d^* \alpha) (x) f(x) + 2 \| \alpha_x \| _{T^* _x M} ^2 f(x) + 2 \mathrm i \langle \alpha_x, \mathrm d _x f \rangle _{T^* _x M \otimes \mathbb C} - \| \alpha_x \| _{T^* _x M} ^2 f(x) = \\
& = \mathrm i (\mathrm d ^* \alpha) (x) \int _{[0,1]} (2 \tau - 1) \, \mathrm d P (\tau) \, f(x) + (\Delta f) (x) \ ,
\end{align*}
whence, finally returning to formula \eqref{Chernoff condition},
\[ \partial_t | _{t=0} [(R_{\alpha, t} f) (x)] = (\Delta f) (x) = (-L^{(U)} f) (x) \]
for every $x \in U$ (the last equality resulting from $f \in \mathcal E \subset \operatorname{Dom} (L^{(U)})$, and in particular from the smoothness of $f$).

Denoting by $A$ the integral $\int _{[0,1]} (2 \tau - 1) \, \mathrm d P (\tau)$ for less visual clutter, the second derivative of $(R_{\alpha, t} ^{(U)} f) (x)$ with respect to $t$ is
\begin{align*}
\partial_s ^2 & (R_{\alpha, s} ^{(U)} f) (x) = \int_U h_\alpha ^{(U)} (s,x,y) \, [\Delta^{(\alpha)} _{(y)}]^2 \, [\chi(x, y) \, \mathrm e ^{-\mathrm i I_P(\alpha) (x,y) - \mathrm i s (\mathrm d^* \alpha) (x) A} f(y)] \, \mathrm d y + \\
& + \int_U h_\alpha ^{(U)} (s,x,y) \, \Delta^{(\alpha)} _{(y)} \, \left[ \chi(x,y) \, \mathrm e ^{-\mathrm i I_P(\alpha) (x,y) - \mathrm i s (\mathrm d^* \alpha) (x) A} (-\mathrm i) (\mathrm d ^* \alpha) (x) \, A \, f(y) \right] \, \mathrm d y + \\
& + \int_U h_\alpha ^{(U)} (s,x,y) \, \chi(x,y) \, \mathrm e ^{-\mathrm i I_P(\alpha) (x,y) - \mathrm i s (\mathrm d^* \alpha) (x) A} (-1) ((\mathrm d ^* \alpha) (x))^2 A^2 f(y) \, \mathrm d y
\end{align*}
whence, using again that $|h_\alpha ^{(U)}| \le h ^{(U)}$ and that $\int_U h ^{(U)} (s,x,y) \, \mathrm d y \le 1$, we obtain the bound
\begin{align*}
| \partial_s ^2 (R_{\alpha, s} ^{(U)} f) (x) | & \le \sup_{y \in \overline U} |[\Delta^{(\alpha)} _{(y)}]^2 \, [\chi(x,y) \, \mathrm e ^{-\mathrm i I_P(\alpha) (x,y)} f(y)]| + \\
& + \sup_{y \in \overline U} |\Delta^{(\alpha)} _{(y)} \, [\chi(x,y) \, \mathrm e ^{-\mathrm i I_P(\alpha) (x,y)} f(y)]| \, |(\mathrm d ^* \alpha) (x)| \int_{[0,1]} (2\tau -1) \, \mathrm d P (\tau) + \\
& + \sup_{y \in \overline U} \chi(x,y) \, |f(y)| \, (\mathrm d ^* \alpha) (x)^2 \left( \int_{[0,1]} (2\tau -1) \, \mathrm d P (\tau) \right)^2 \ .
\end{align*}

From the compactness of $\overline U$ and the continuity (and in fact smoothness) of all the functions in the right hand side, we get that each of the three terms that make up the latter is bounded with respect to $x \in \overline U$. We thus obtain that there exists a constant $C>0$ (depending of all the objects involved, of course, i.e. $\alpha$, $\chi$, $P$ and $U$) such that $| \partial_s ^2 (R_s ^{(U)} f) (x) | \le C$, whence, returning to formula (\ref{Taylor expansion}) with all the results obtained above,
\[ \| R_{\alpha, t} ^{(U)} f - f - (-L^{(U)} f) \, t \| _{C(\overline U)} \le C \frac {t^2} 2 \ , \]
hence the last hypothesis in Chernoff's theorem that had to be checked is now immediate (on $\partial U$ the functions $f$, $L^{(U)} f$ and $R_{\alpha, t} ^{(U)} f$ vanish, because $f \in \mathcal E$, so the behaviour of the functions on the boundary does not alter the conclusion).

Since this was the last thing to check, we deduce that we may apply Chernoff's theorem, thus obtaining that $\mathrm e^{-tL^{(U)}} = \lim _{k \to \infty} \big( R_{\alpha, \frac t k} ^{(U)} \big) ^k$ strongly in $C(\overline U)$, which is exactly what we were trying to prove.
\end{proof}

It becomes clear now where the compactness of $\overline U$ mattered: at the very end of the proof, where we have used that continuous functions are bounded on compact subsets; this clarifies why we were led to work on relatively compact subsets and not directly on the whole manifold.

For every $k \in \mathbb N$, let us now define the "approximations" $S_{P,t,k} (\alpha) : \mathcal C_t \to \mathbb R$ by
\begin{align*}
S_{P,t,k} (\alpha) (c) & = \sum _{j=0} ^{2^k - 1} I_P (\alpha) \left( c \left( \frac {jt} {2^k} \right), c \left( \frac {(j+1)t} {2^k} \right) \right) + \\
& + \frac t {2^k} (\mathrm d ^* \alpha) \left( c \left( \frac {jt} {2^k} \right) \right) \int _{[0,1]} (2 \tau - 1) \, \mathrm d P (\tau) \ .
\end{align*}

So far, $\rho_{\alpha, t} ^{(U)}$ has been obtained by a very abstract procedure (section \ref{rho}), which makes its use in concrete calculations and the study of its properties very difficult. The following theorem remedies this situation, providing us with a concrete understanding of $\rho_{\alpha, t} ^{(U)}$ as the limit of a sequence of functions given by explicit formulae.

\begin{theorem} \label{convergence on regular domains}
$\lim _{k \to \infty} \mathrm e ^{\mathrm i S_{P,t,k} (\alpha)} \big| _{\mathcal C_t (\overline U)} = \rho_{\alpha, t} ^{(U)}$ in $L^2(\mathcal C_t(\overline U), w_t ^{(U)})$, uniformly with respect to $t$ in bounded subsets of $(0, \infty)$, and uniformly with respect to $x_0 \in U$.
\end{theorem}
\begin{proof}
We shall reduce the problem to the application of theorem \ref{the technical theorem}; for notational simplicity, we shall write $\mathrm e ^{\mathrm i S_{P,t,k} (\alpha)}$ instead of $\mathrm e ^{\mathrm i S_{P,t,k} (\alpha)} \big| _{\mathcal C_t (\overline U)}$. Clearly,
\begin{gather} \label{limit in L^2}
\nonumber \left\| \mathrm e ^{\mathrm i S_{P,t,k} (\alpha)} - \rho_{\alpha,t} ^{(U)} \right\| _{L^2(\mathcal C_t (\overline U))} ^2 = \\
= \left\| \mathrm e ^{\mathrm i S_{P,t,k} (\alpha)} \right\| _{L^2(\mathcal C_t (\overline U))} ^2 - \left\langle \mathrm e ^{\mathrm i S_{P,t,k} (\alpha)}, \rho_{\alpha,t} ^{(U)} \right\rangle _{L^2(\mathcal C_t (\overline U))} - \left\langle \rho_{\alpha,t} ^{(U)}, \mathrm e ^{\mathrm i S_{P,t,k} (\alpha)} \right\rangle _{L^2(\mathcal C_t (\overline U))} + \left\| \rho_{\alpha,t} ^{(U)} \right\| _{L^2(\mathcal C_t (\overline U))} ^2 \ ,
\end{gather}
and the first term of the sum is obviously $w_t ^{(U)} (\mathcal C_t (\overline U))$.

In order to evaluate the third term (and thus the second, by conjugation), let us define the function $\chi_k : \mathcal C_t (\overline U) \to [0,1]$ by
\[ \chi_k (c) = \chi \left( c(0), c\left( \frac t {2^k} \right) \right) \chi \left( c \left( \frac t {2^k} \right), c \left( \frac {2t} {2^k} \right) \right) \dots \chi \left( c \left( \frac {(2^k - 1)t} {2^k} \right), c(t) \right) \]
for every $k \ge 0$. This allows us to write that
\begin{gather*}
\lim_{k \to \infty} \left\langle \rho_{\alpha,t} ^{(U)}, \mathrm e ^{\mathrm i S_{P,t,k} (\alpha)} \right\rangle _{L^2(\mathcal C_t (\overline U))} = \lim_{k \to \infty} \int _{\mathcal C_t (\overline U)} \rho_{\alpha,t} ^{(U)} \, \mathrm e ^{- \mathrm i S_{P,t,k} (\alpha)} \, \mathrm d w_t ^{(U)} = \\
= \lim_{k \to \infty} W_{\alpha, t} ^{(U)} (\mathrm e ^{- \mathrm i S_{P,t,k} (\alpha)}) = \lim_{k \to \infty} W_{\alpha, t} ^{(U)} (\chi_k \, \mathrm e ^{- \mathrm i S_{P,t,k} (\alpha)}) + \lim_{k \to \infty} W_{\alpha, t} ^{(U)} ((1-\chi_k) \, \mathrm e ^{- \mathrm i S_{P,t,k} (\alpha)}) \ .
\end{gather*}

Using theorem \ref{the technical theorem}, the first term becomes
\begin{gather*}
\lim _{k \to \infty} \int_U \mathrm d x_1 \, h_\alpha ^{(U)} \left( \frac t {2^k}, x_0, x_1 \right) \, \chi(x_0, x_1) \, \mathrm e ^{- \mathrm i I_P (\alpha) (x_0, x_1) - \frac t {2^k} \mathrm i (\mathrm d ^* \alpha) (x_0) \, [2 M_1(P) - 1]} \dots \\
\dots \int_U \mathrm d x_{2^k} \, h_\alpha ^{(U)} \left( \frac t {2^k}, x_{2^k - 1}, x_{2^k} \right) \, \chi(x_{2^k-1}, x_{2^k}) \, \mathrm e ^{- \mathrm i I_P (\alpha) (x_{2^k - 1}, x_{2^k}) - \frac t {2^k} \mathrm i (\mathrm d ^* \alpha) (x_{2^k - 1}) \, [2 M_1(P) - 1]} = \\
= \left[ \lim _{k \to \infty} \left( R_{\alpha, \frac t {2^k}} ^{(U)} \right) ^{2^k} 1 \right] (x_0) = (\mathrm e ^{-t L^{(U)}} \, 1) (x_0) = \int_U h ^{(U)} (t, x_0, x) \, \mathrm d x = w_t ^{(U)} (\mathcal C_t (\overline U)) \ .
\end{gather*}

Notice that the limit that we have just obtained is uniform with respect to $t$ in bounded subsets of $(0, \infty)$, as a consequence of the uniformity obtained in theorem \ref{the technical theorem}.

Since the diamagnetic inequality $| h_\alpha ^{(U)} | \le h ^{(U)}$ implies that $| W_{\alpha, t} ^{(U)} (f) | \le \int _{\mathcal C_t (\overline U)} |f| \, \mathrm d w_t ^{(U)}$ for every $f \in L^1 (\mathcal C_t (\overline U))$, we obtain for the second term that
\begin{align*}
0 \le \lim_{k \to \infty} | W_{\alpha, t} ^{(U)} ((1-\chi_k) \, \mathrm e ^{- \mathrm i S_{P,t,k} (\alpha)}) | \le \lim_{k \to \infty} \int _{\mathcal C_t (\overline U)} 1 - \chi_k \, \mathrm d w_t ^{(U)} = \\
= w_t ^{(U)} (\mathcal C_t (\overline U)) - \lim_{k \to \infty} W_{0, t} ^{(U)} (\chi_k \, \mathrm e ^{-\mathrm i S_{P,t,k} (0)}) = 0 \ ,
\end{align*}
where the last limit has been obtained with the same reasoning as above performed for $\alpha=0$. Notice that this limit, too, is uniform with respect to $t$ in bounded subsets of $(0, \infty)$, for the same reason as above.

So far, then, formula \eqref{limit in L^2} has given us
\begin{align*}
0 \le \limsup _{k \to \infty} \left\| \mathrm e ^{\mathrm i S_{P,t,k} (\alpha)} - \rho_{\alpha,t} ^{(U)} \right\| _{L^2(\mathcal C_t (\overline U))} ^2 = \left\| \rho_{\alpha,t} ^{(U)} \right\| _{L^2(\mathcal C_t (\overline U))} ^2 - w_t ^{(U)} (\mathcal C_t (\overline U)) \le \\
\le w_t ^{(U)} (\mathcal C_t (\overline U)) \, \| \rho_{\alpha,t} ^{(U)} \| _{L^\infty (\mathcal C_t (\overline U))} ^2 - w_t ^{(U)} (\mathcal C_t (\overline U)) \ ,
\end{align*}
the limit superior being uniform with respect to $t$ in bounded subsets of $(0, \infty)$. But
\begin{align*}
\| \rho_{\alpha,t} ^{(U)} \| _{L^\infty (\mathcal C_t (\overline U))} & = \sup _{\| f \| _{L^1 (\mathcal C_t (\overline U))} \le 1} \left| \int _{\mathcal C_t (\overline{U})} \rho_{\alpha, t} ^{(U)} \, f \, \mathrm d w_t ^{(U)} \right| = \sup _{\| f \| _{L^1 (\mathcal C_t (\overline U))} \le 1} | W_{\alpha, t} ^{(U)} (f) | \le \\
& \le \sup _{\| f \| _{L^1 (\mathcal C_t (\overline U))} \le 1} \| f \| _{L^1 (\mathcal C_t (\overline U))} \le 1 \ ,
\end{align*}
whence we conclude that $0 \le \limsup _{k \to \infty} \left\| \mathrm e ^{\mathrm i S_{P,t,k} (\alpha)} - \rho_{\alpha, t} ^{(U)} \right\| _{L^2(\mathcal C_t (\overline U))} ^2 \le 0$ and, since
\[ 0 \le \liminf _{k \to \infty} \left\| \mathrm e ^{\mathrm i S_{P,t,k} (\alpha)} - \rho_{\alpha, t} ^{(U)} \right\| _{L^2(\mathcal C_t (\overline U))} ^2 \le \limsup _{k \to \infty} \left\| \mathrm e ^{\mathrm i S_{P,t,k} (\alpha)} - \rho_{\alpha, t} ^{(U)} \right\| _{L^2(\mathcal C_t (\overline U))} ^2 \le 0 \ , \]
we conclude that, indeed, $\lim _{k \to \infty} \mathrm e ^{\mathrm i S_{P,t,k} (\alpha)} = \rho_{\alpha, t} ^{(U)}$ in $L^2 (\mathcal C_t (\overline U))$, the limit being uniform with respect to $t$ in bounded subsets of $(0, \infty)$, as desired.
\end{proof}

\section{A unitary group and its generator} \label{the Stratonovich integral}

Let us now consider an exhaustion $M = \bigcup _{j \in \mathbb N} U_j$ of $M$ with regular domains (it exists as a consequence of proposition 2.28 in \cite{Lee13}). For notational simplicity, let us write $\rho_{\alpha, t} ^{(j)}$ instead of $\rho_{\alpha, t} ^{(U_j)}$, $h_\alpha ^{(j)}$ instead of $h_\alpha ^{(U_j)}$, and $w_t ^{(j)}$ instead of $w_t ^{(U_j)}$. So far we know that $\mathrm e ^{\mathrm i S_{P,t,k}} \big| _{\mathcal C_t (\overline {U_j})} \to \rho_{\alpha, t} ^{(j)}$ in $L^2 (\mathcal C_t (\overline {U_j}), w_t ^{(j)})$ for every $j \in \mathbb N$.

\begin{lemma}
The subset $\mathcal C_t (\overline {U_j})$ is closed in $\mathcal C_t$ for every $j \ge 0$. Similarly, $\mathcal C_t (\overline {U_i})$ is closed in $\mathcal C_t (\overline {U_j})$ for every $i \le j$.
\end{lemma}
\begin{proof}
We shall prove only the first statement, the proof of the second being similar. The evaluation map $\operatorname{ev} : [0,t] \times \mathcal C_t \to M$ defined by $\operatorname{ev} (s, \gamma) = \gamma(s)$ is obviously continuous, whence
\[ \mathcal C_t (\overline {U_j}) = \{ \gamma \in C_t \mid \gamma(s) \in \overline {U_j} \ \forall s \in [0,t] \} = \bigcap _{s \in [0,t]} \operatorname{ev} (s, \cdot) ^{-1} (\overline {U_j}) \]
is obviously closed.
\end{proof}

The following lemma is as important as it is trivial.

\begin{lemma}
If $i \le j$ then $\rho_{\alpha, t} ^{(j)} \big| _{\mathcal C_t (\overline {U_i})} = \rho_{\alpha, t} ^{(i)}$ almost everywhere on $\mathcal C_t (\overline {U_i})$ with respect to the Wiener measure $w_t ^{(i)}$.
\end{lemma}
\begin{proof}
It is a consequence of the inequality $h ^{(U_i)} \le h ^{(U_j)} | _{U_i}$ (which, in turn, is a consequence of the parabolic minimum principle applied to the heat equation on $U_i$) that $w_t ^{(i)} \le w_t ^{(j)} \big| _{\mathcal C_t (\overline {U_i})}$ (for details see \cite{BP11}). It then follows that for every $k \in \mathbb N$ we have
\begin{align*}
& \| \rho_{\alpha, t,} ^{(j)} | _{\mathcal C_t (\overline{U_i})} - \rho_{\alpha, t} ^{(i)} \| _{L^2 (\mathcal C_t (\overline{U_i}))} \le \\
& \le \| \rho_{\alpha, t} ^{(j)} | _{\mathcal C_t (\overline{U_i})} - \mathrm e ^{\mathrm i S_{P,t,k}} | _{\mathcal C_t (\overline{U_i})} \| _{L^2 (\mathcal C_t (\overline{U_i}))} + \| \mathrm e ^{\mathrm i S_{P,t,k}} | _{\mathcal C_t (\overline{U_i})} - \rho_{\alpha, t} ^{(i)} \| _{L^2 (\mathcal C_t (\overline{U_i}))} = \\
& = \left( \int _{\mathcal C_t (\overline{U_i})} | \rho_{\alpha, t} ^{(j)} (c) - \mathrm e ^{\mathrm i S_{P,t,k}} (c) | ^2 \, \mathrm d w_t ^{(i)} (c) \right) ^{\frac 1 2} + \| \mathrm e ^{\mathrm i S_{P,t,k}} (c) | _{\mathcal C_t (\overline{U_i})} - \rho_{\alpha, t} ^{(i)} \| _{L^2 (\mathcal C_t (\overline{U_i}))} \le \\
& \le \left( \int _{\mathcal C_t (\overline{U_i})} | \rho_{\alpha, t} ^{(j)} (c) - \mathrm e ^{\mathrm i S_{P,t,k}} (c) | ^2 \, \mathrm d w_t ^{(j)} | _{\mathcal C_t (\overline{U_i})} (c) \right) ^{\frac 1 2} + \| \mathrm e ^{\mathrm i S_{P,t,k}} | _{\mathcal C_t (\overline{U_i})} - \rho_{\alpha, t} ^{(i)} \| _{L^2 (\mathcal C_t (\overline{U_i}))} \le \\
& \le \left( \int _{\mathcal C_t (\overline{U_j})} | \rho_{\alpha, t} ^{(j)} (c) - \mathrm e ^{\mathrm i S_{P,t,k}} (c) | ^2 \, \mathrm d w_t ^{(j)} (c) \right) ^{\frac 1 2} + \| \mathrm e ^{\mathrm i S_{P,t,k}} | _{\mathcal C_t (\overline{U_i})} - \rho_{\alpha, t} ^{(i)} \| _{L^2 (\mathcal C_t (\overline{U_i}))} = \\
& = \| \rho_{\alpha, t} ^{(j)} - \mathrm e ^{\mathrm i S_{P,t,k}} | _{\mathcal C_t (\overline{U_j})} \| _{L^2 (\mathcal C_t (\overline{U_j}))} + \| \mathrm e ^{\mathrm i S_{P,t,k}} | _{\mathcal C_t (\overline{U_i})} - \rho_{\alpha, t} ^{(i)} \| _{L^2 (\mathcal C_t (\overline{U_i}))} \ ,
\end{align*}
whence the conclusion is clear using theorem \ref{convergence on regular domains}.
\end{proof}

The equality $\rho_{\alpha, t} ^{(j)} \big| _{\mathcal C_t (\overline {U_i})} = \rho_{\alpha, t} ^{(i)}$ almost everywhere for every $i \le j$ implies the existence of the pointwise limit $\lim _{j \to \infty} \rho_{\alpha, t} ^{(j)}$, which we shall denote by $\rho _{\alpha, t}$. It will be a measurable function (as the pointwise limit of a sequence of measurable functions), and it will be bounded by $1$ almost everywhere, because all the functions in the sequence are so. Therefore, it will be an element of $L^\infty (\mathcal C_t, w_t)$. Using the argument in the above lemma, one may show that $\rho_{\alpha, t} \big| _{\mathcal C_t (\overline {U_j})} = \rho_{\alpha, t} ^{(j)}$ for every $j \ge 0$, as elements from $L^\infty (\mathcal C_t, w_t ^{(j)})$.

After all these preliminary results, we may finally prove one of the core results of this work.

\begin{theorem}
\[ \lim _{k \to \infty} \| \mathrm e ^{\mathrm i S_{P,t,k} (\alpha)} - \rho_{\alpha, t} \| _{L^2 (\mathcal C_t)} = 0 \]
uniformly with respect to $t \in (0,T]$, for every $T>0$.
\end{theorem}
\begin{proof}
Let $\varepsilon > 0$. Using the fact that $\rho_{\alpha, t} \big| _{\mathcal C_t (\overline {U_j})} = \rho_{\alpha, t} ^{(j)}$ for every $j \ge 0$, we may write that
\begin{align*}
\| \mathrm e ^{\mathrm i S_{P,t,k} (\alpha)} - \rho_{\alpha, t} \| _{L^2 (\mathcal C_t, w_t)} ^2 & = \| \mathrm e ^{\mathrm i S_{P,t,k} (\alpha)} - \rho_{\alpha, t} \| _{L^2 (\mathcal C_t (\overline {U_j}), w_t)} ^2 + \| \mathrm e ^{\mathrm i S_{P,t,k} (\alpha)} - \rho_{\alpha, t} \| _{L^2 (\mathcal C_t \setminus \mathcal C_t (\overline {U_j}), w_t)} ^2 \le \\
& \le \| \mathrm e ^{\mathrm i S_{P,t,k} (\alpha)} - \rho_{\alpha, t} ^{(j)} \| _{L^2 (\mathcal C_t (\overline {U_j}), w_t)} ^2 + 4 w_t(\mathcal C_t \setminus \mathcal C_t (\overline {U_j})) = \\
& \le \| \mathrm e ^{\mathrm i S_{P,t,k} (\alpha)} - \rho_{\alpha, t} ^{(j)} \| _{L^2 (\mathcal C_t (\overline {U_j}), w_t)} ^2 + 4 [w_t(\mathcal C_t) - w_t ^{(j)} (\mathcal C_t (\overline {U_j}))] \ .
\end{align*}

The first summand, in turn, may be written as
\begin{align*}
\| \mathrm e ^{\mathrm i S_{P,t,k} (\alpha)} & - \rho_{\alpha, t} ^{(j)} \| _{L^2 (\mathcal C_t (\overline {U_j}), w_t)} ^2 = \\
& = \int _{\mathcal C_t (\overline {U_j})} | \mathrm e ^{\mathrm i S_{P,t,k} (\alpha)} - \rho_{\alpha, t} ^{(j)} | ^2 \, \mathrm d w_t ^{(j)} + \int _{\mathcal C_t (\overline {U_j})} | \mathrm e ^{\mathrm i S_{P,t,k} (\alpha)} - \rho_{\alpha, t} ^{(j)} | ^2 \, \mathrm d (w_t - w_t ^{(j)}) \le \\
& \le \| \mathrm e ^{\mathrm i S_{P,t,k} (\alpha)} - \rho_{\alpha, t} ^{(j)} \| _{L^2 (\mathcal C_t (\overline {U_j}), w_t ^{(j)})} ^2 + 4 [w_t (\mathcal C_t (\overline {U_j})) - w_t ^{(j)} (\mathcal C_t (\overline {U_j}))] \le \\
& \le \| \mathrm e ^{\mathrm i S_{P,t,k} (\alpha)} - \rho_{\alpha, t} ^{(j)} \| _{L^2 (\mathcal C_t (\overline {U_j}), w_t ^{(j)})} ^2 + 4 [w_t (\mathcal C_t) - w_t ^{(j)} (\mathcal C_t (\overline {U_j}))] \ ,
\end{align*}
therefore we obtain
\begin{align*}
\| \mathrm e ^{\mathrm i S_{P,t,k} (\alpha)} - \rho_{\alpha, t} \| _{L^2 (\mathcal C_t, w_t)} ^2 & \le \| \mathrm e ^{\mathrm i S_{P,t,k} (\alpha)} - \rho_{\alpha, t} ^{(j)} \| _{L^2 (\mathcal C_t (\overline {U_j}), w_t ^{(j)})} ^2 + 8 [w_t (\mathcal C_t) - w_t ^{(j)} (\mathcal C_t (\overline {U_j}))] \ .
\end{align*}

We shall see next how to choose $j$ and $k$ such that the right hand side should be less than $\varepsilon$, independently of $t \in (0, T]$.

The second term contains
\[ w_t (\mathcal C_t) - w_t ^{(j)} (\mathcal C_t (\overline {U_j})) = \int _M h(t, x_0, x) \, \mathrm d x - \int _{U_j} h^{(j)} (t, x_0, x) \, \mathrm d x \ . \]
Since $h^{(j)} \to h$ pointwise (and monotonically), there exists a $j_{\varepsilon, t}$ for every $t > 0$ such that
\[ \left| \int _M h(t, x_0, x) \, \mathrm d x - \int _{U_j} h^{(j)} (t, x_0, x) \, \mathrm d x \right| = \left| \langle \delta_{x_0}, \mathrm e ^{-t L} 1 \rangle - \langle \delta_{x_0}, \mathrm e ^{-t L^{(j)}} 1 \big| _{\overline {U_j}} \rangle \right| < \frac \varepsilon {16} \]
for every $j \ge j_{\varepsilon, t}$, where $\langle \cdot, - \rangle$ denotes the dual pairing between $C_b (M)$ (and, respectively, $C(\overline {U_j})$) and its topological dual. Since those heat semigroups are strongly continuous, the expression inside the absolute value is continuous with respect to $t \in [0,T]$, therefore there exists some small open neighbourhood $V_{\varepsilon, t} \subseteq [0,T]$ of $t$ in $[0,T]$ such that
\[ \left| (\mathrm e ^{-s L} 1) (x_0) - (\mathrm e ^{-s L^{(j)}} 1 \big| _{\overline {U_j}}) (x_0) \right| < \frac \varepsilon {16} \]
for every $j \ge j_{\varepsilon, t}$, uniformly with respect to $s \in V_{\varepsilon, t}$. Since $[0,T]$ is compact, it may be written as a finite union $[0,T] = \bigcup _{i=1} ^{N_{\varepsilon}} V_{\varepsilon, t_i}$. We shall then choose $j = \max \{ j_{\varepsilon, t_1}, \dots, j_{\varepsilon, t_N} \}$ and we shall thus have
\[ \left| \int _M h(t, x_0, x) \, \mathrm d x - \int _{U_j} h^{(j)} (t, x_0, x) \, \mathrm d x \right| < \frac \varepsilon {16} \]
for every $t \in (0,T]$.

With $j$ so chosen, we shall use theorem \ref{convergence on regular domains} on $\mathcal C_t (\overline {U_j})$ in order to find $k_\varepsilon$ such that
\[ \| \mathrm e ^{\mathrm i S_{P,t,k} (\alpha)} - \rho_{\alpha, t} ^{(j)} \| _{L^2 (\mathcal C_t (\overline {U_j}), w_t ^{(j)})} ^2 < \frac \varepsilon 2 \]
for every $k \ge k_\varepsilon$, uniformly with respect to $t \in (0,T]$.

Combining these two majorizations, we obtain that
\[ \| \mathrm e ^{\mathrm i S_{P,t,k} (\alpha)} - \rho_{\alpha, t} \| _{L^2 (\mathcal C_t)} ^2 < \varepsilon \]
for every $k \ge k_\varepsilon$, uniformly with respect to $t \in (0,T]$, whence the conclusion is clear.
\end{proof}

\begin{corollary}
$\rho_{\alpha, t}$ does not depend on the exhaustion with regular domains used.
\end{corollary}
\begin{proof}
Clear, because the approximations $\mathrm e ^{\mathrm i S_{P,t,k} (\alpha)}$ do not depend on any exhaustion.
\end{proof}

Let us break the pace of our exposition and make a comment regarding the strategy chosen to define $\rho_{\alpha, t}$: we could have defined it directly on $\mathcal C_t$, as a representative from $L^\infty (\mathcal C_t)$ of the continuous linear functional $W_{\alpha, t}$ on $L^1 (\mathcal C_t)$, exactly like in section \ref{rho}, thus skipping the long intermediate and purely ancillary step involving the exhaustion of $M$. The problem is that in this approach the connection with the sequence of approximating exponentials becomes much more difficult to exhibit; one may show that both approaches lead to the same $\rho_{\alpha, t}$, but the proof is tedious and its usefulness not clear, therefore we shall not bother with it anymore.

We have obtained that $\rho_{\alpha, t}$ is the limit of a sequence of exponentials with imaginary exponents. It is reasonable to ask whether $\rho_{\alpha, t}$ itself has such a form and, if the answer is affirmative, to study its exponent. The answer to this question (and the moral justification of all the effort spent in obtaining all the technical results so far) is given by theorem \ref{existence of the Stratonovich integral}. In order to state it, though, we need the following lemma, the proof of which we omit in order to save space, since it is fairly elementary and found in undergraduate-level textbooks.

\begin{lemma}
Let $(X, \mu)$ be a space endowed with a $\sigma$-finite measure, and let $T : L^2 (X) \to L^2 (X)$ be a bounded linear operator. If $T(\varphi f) = \varphi T(f)$ for every $\varphi \in L^\infty (X)$ and $f \in L^2 (X)$, then there exists a function $\theta \in L^\infty (X)$ such that $T(f) = \theta f$ for every $f \in L^2 (X)$.
\end{lemma}

\begin{theorem} \label{existence of the Stratonovich integral}
There exists a unique real-valued function $\operatorname{Strat} _t (\alpha) \in L^0 (\mathcal C_t)$ such that $\rho_{\alpha, t} = \mathrm e ^{\mathrm i \operatorname{Strat} _t (\alpha)}$.
\end{theorem}
\begin{proof}
We have already shown that $\lim _{k \to \infty} \mathrm e ^{\mathrm i S_{P,t,k} (\alpha)} = \rho_{\alpha, t}$ in $L^2(\mathcal C_t)$, therefore there exists a sub-sequence of functions $\big( \mathrm e ^{\mathrm i S_{P,k_j, t} (\alpha)} \big) _{j \ge 0}$ such that $\lim_{j \to \infty} \mathrm e ^{\mathrm i S_{P,k_j, t} (\alpha)} = \rho_{\alpha, t}$ pointwise almost everywhere, whence in particular $| \rho_{\alpha, t} | = 1$ almost everywhere.

Let us consider the multiplication operators $U_s : L^2(\mathcal C_t) \to L^0(\mathcal C_t)$ given by $U_s (f) = \rho_{s \alpha, t} f$ for every $s \in \mathbb R$. Since $| \rho_{s \alpha, t} | = 1$, we have that $U_s (L^2(\mathcal C_t)) \subseteq L^2(\mathcal C_t)$ and $U_s$ is unitary.

Let us show that these operators form a group. On the one hand, it is obvious that $S_{P,t,k}((s+s') \alpha) = S_{P,t,k}(s\alpha) + S_{P,t,k}(s' \alpha)$ for every $k \ge 0$ and $s,s' \in \mathbb R$. On the other hand, it is clear that, repeatedly passing to sub-sequences as many times as necessary, there exists a sequence of numbers $(k_j) _{j \ge 0} \subseteq \mathbb N$ such that $\mathrm e ^{\mathrm i S_{P,k_j,t} ((s+s') \alpha)} \to \rho_{(s+s') \alpha, t}$, $\mathrm e ^{\mathrm i S_{P,k_j,t} (s \alpha)} \to \rho_{s \alpha, t}$ and $\mathrm e ^{\mathrm i S_{P,k_j,t} (s' \alpha)} \to \rho_{s' \alpha, t}$ for $j \to \infty$, whence it follows that
\begin{align*}
U_{s+s'} f & = \rho_{(s+s') \alpha, t} f = \lim _{j \to \infty} \mathrm e ^{\mathrm i S_{P,k_j,t} ((s+s') \alpha)} f = \lim _{j \to \infty} \mathrm e ^{\mathrm i S_{P,k_j,t} (s \alpha)} \lim _{j \to \infty} \mathrm e ^{\mathrm i S_{P,k_j,t} (s' \alpha)} f = \\
& = \rho_{s \alpha, t} \rho_{s' \alpha, t} f = U_s U_{s'} f
\end{align*}
for every $f \in L^2(\mathcal C)$ and $s,s' \in \mathbb R$.

In order to show that this group is strongly continuous, let us remember that $\mathcal C_t$ is separable (see \cite{Michael61}), therefore $L^2 (\mathcal C_t)$ is separable too, so that strong continuity is equivalent to weak measurability. If $f,g \in L^2 (\mathcal C_t)$, then
\[\langle U_s f, g \rangle _{L^2(\mathcal C_t)} = \langle \rho_{s \alpha, t} f, g \rangle _{L^2(\mathcal C_t)} = \lim _{k \to \infty} \langle \mathrm e ^{\mathrm i S_{P,t,k} (s \alpha)} f, g \rangle _{L^2(\mathcal C_t)} \ ,\]
which is the pointwise limit of a sequence of functions that obviously depend measurably on $s$, therefore it is measurable itself.

Denoting by $\operatorname{Mul} (L^2(\mathcal C_t))$ the algebra of multiplication operators on $L^2(\mathcal C_t)$, we have obtained that the map $\mathbb R \ni s \mapsto U_s \in \operatorname{Mul} (L^2(\mathcal C_t))$ is a strongly continuous unitary one-parameter group, therefore we deduce (by Stone's theorem) that there exists a unique self-adjoint operator $\operatorname{Strat} _t (\alpha)$ such that $\rho_{s \alpha, t} = \mathrm e ^{\mathrm i s \operatorname{Strat} _t (\alpha)}$ for every $s \in \mathbb R$. Let us show that $\operatorname{Strat} _t (\alpha)$ is a multiplication operator itself.

Since $\operatorname{Strat} _t (\alpha)$ is self-adjoint, we may consider its resolvent operator
\[ R_{- \mathrm i} (\operatorname{Strat} _t (\alpha)) = [\operatorname{Strat} _t (\alpha) - (-\mathrm i)] ^{-1} \]
at $- \mathrm i$, given by Laplace's formula
\[ R_{- \mathrm i} (\operatorname{Strat} _t (\alpha)) f = -\mathrm i \int _0 ^\infty \mathrm e ^{-s} \, U_s f \, \mathrm d s \]
for every $f \in L^2 (\mathcal C_t)$. The operator $R_{- \mathrm i} (\operatorname{Strat} _t (\alpha))$ is bounded, and an elementary application of Fubini's theorem shows that
\[ R_{- \mathrm i} (\operatorname{Strat} _t (\alpha)) \operatorname{Mul} (\varphi) = \operatorname{Mul} (\varphi) R_{- \mathrm i} (\operatorname{Strat} _t (\alpha)) \ , \]
where $\operatorname{Mul} (\varphi)$ is the multiplication operator by $\varphi \in L^\infty (\mathcal C_t)$. Since $\varphi$ is arbitrary, it follows that $R_{- \mathrm i} (\operatorname{Strat} _t (\alpha))$ is given by the multiplication by a unique essentially bounded function which we shall denote by $s_t (\alpha)$. This function is non-zero almost everywhere: if $E \subseteq \mathcal C_t$ is such that $w_t (E) \ne 0$ and $s_t | _{E} = 0$, then the characteristic function $1_E$ is an eigenvector for $R_{- \mathrm i} (\operatorname{Strat} _t (\alpha))$ with eigenvalue $0$, whence
\[ (\mathrm i + S(\alpha)) \, R_{- \mathrm i} (\operatorname{Strat} _t (\alpha)) 1_E = 0 \ , \]
which contradicts invertibility.

Since
\[ R_{- \mathrm i} (\operatorname{Strat} _t (\alpha)) = [\mathrm i + \operatorname{Strat} _t (\alpha)] ^{-1} \ , \]
it follows that $\operatorname{Strat} _t (\alpha) = \mathrm i + \frac 1 {s_t (\alpha)}$, hence the operator $\operatorname{Strat} _t (\alpha)$ is given by the multiplication by a unique element of $L^0 (\mathcal C_t)$, that we shall keep denoting by $\operatorname{Strat} _t (\alpha)$.

The fact that the function $\operatorname{Strat} _t (\alpha)$ is real-valued follows from the fact that the operator $\operatorname{Strat} _t (\alpha)$ is self-adjoint.
\end{proof}

When we constructed the functions $S_{P,t,k} (\alpha)$, we did it in order for the functions $\mathrm e ^{\mathrm i S_{P,t,k} (\alpha)}$ to approximate $\rho_{\alpha, t} = \mathrm e ^{\mathrm i \operatorname{Strat} _t (\alpha)}$ in $L^2 (\mathcal C_t)$. We shall see now that this approximation property extends, even though in a weaker form, to the exponents.

\begin{theorem} \label{approximation of the Stratonovich integral}
$\lim _{k \to \infty} S_{P,t,k} (\alpha) = \operatorname{Strat} _t (\alpha)$ in measure, uniformly with respect to $t$ in bounded subsets of $(0, \infty)$.
\end{theorem}
\begin{proof}
Using the notations and the integral formula for the resolvent operator from the previous theorem, and taking into account the fact that $1 \in L^2 (\mathcal C_t)$, and that all the operators appearing in this proof are multiplication operators, we have
\begin{gather*}
\| R_{- \mathrm i} (\operatorname{Strat} _t (\alpha)) - R_{- \mathrm i} (S_{P,t,k} (\alpha)) \| _{L^2(\mathcal C_t)} = \| R_{- \mathrm i} (\operatorname{Strat} _t (\alpha)) \, 1 - R_{- \mathrm i} (S_{P,t,k} (\alpha)) \, 1 \| _{L^2(\mathcal C_t)} \le \\
\le \int _0 ^\infty \mathrm e ^{-s} \, \| \rho_{s \alpha, t} \, 1 - \mathrm e ^{\mathrm i S_{P,t,k} (s \alpha)} \, 1 \| _{L^2(\mathcal C_t)} \, \mathrm d s = \int _0 ^\infty \mathrm e ^{-s} \, \| \rho_{s \alpha, t} - \mathrm e ^{\mathrm i S_{P,t,k} (s \alpha)} \| _{L^2(\mathcal C_t)} \, \mathrm d s \ .
\end{gather*}
Since
\[ \| \rho_{s \alpha, t} - \mathrm e ^{\mathrm i S_{P,t,k} (s \alpha)} \| _{L^2(\mathcal C_t)} \le \| \rho_{s \alpha, t} - \mathrm e ^{\mathrm i S_{P,t,k} (s \alpha)} \| _{L^\infty(\mathcal C_t)} \sqrt {w_t (\mathcal C_t)} \le 2 \sqrt {w_t (\mathcal C_t)} \ , \]
using the dominated convergence theorem and the fact shown above that $\| \rho_{s \alpha, t} - \mathrm e ^{\mathrm i S_{P,t,k} (s \alpha)} \| _{L^2(\mathcal C_t)} \to 0$ uniformly with respect to $t$ in bounded subsets of $(0, \infty)$, we obtain that $\| R_{- \mathrm i} (\operatorname{Strat} _t (\alpha)) - R_{- \mathrm i} (S_{P,t,k} (\alpha)) \| _{L^2(\mathcal C_t)} \to 0$ uniformly with respect to $t$ in bounded subsets of $(0, \infty)$, and therefore that $\frac 1 {\mathrm i + S_{P,t,k} (\alpha)} \to \frac 1 {\mathrm i + \operatorname{Strat} _t (\alpha)}$ in $L^2 (\mathcal C_t)$, and in particular in measure, when $k \to \infty$, uniformly with respect to $t$ in bounded subsets of $(0, \infty)$, whence the conclusion is immediate.
\end{proof}

We shall see in detail, in the next section, that $\operatorname{Strat} _t$ is the Stratonovich stochastic integral. The fact that it is the limit in measure of the sequence of approximations used above was already known; what is new is that it stems into existence as the generator of the unitary group considered above (or, giving up rigour, it is the "logarithm" of the function $\rho_{\alpha, t}$). This latter fact suggests that $\rho_{\alpha, t}$, being the imaginary exponential of a stochastic line integral, may be viewed as a sort of parallel transport - namely the stochastic parallel transport in the trivial vector bundle $M \times \mathbb C$. These considerations will form the object of a separate study in a more geometric setting, which is not included in the present work.

\begin{corollary}
The map $\Omega^1 (M) \ni \alpha \mapsto \operatorname{Strat} _t (\alpha) \in L^0 (\mathcal C_t)$ is $\mathbb R$-linear.
\end{corollary}
\begin{proof}
\[\operatorname{Strat} _t (\alpha + \beta) = \lim _{k \to \infty} S_{P,t,k} (\alpha + \beta) = \lim _{k \to \infty} S_{P,t,k} (\alpha) + \lim _{k \to \infty} S_{P,t,k} (\beta) = \operatorname{Strat} _t (\alpha) + \operatorname{Strat} _t (\beta) \ ,\]
where all the limits are considered in measure. A similar argument shows that $\operatorname{Strat} _t (r \alpha) = r \operatorname{Strat} _t (\alpha)$, for every $r \in \mathbb R$.
\end{proof}

Let us end this section with a remark regarding the strategy chosen to construct $\operatorname{Strat} _t (\alpha)$: since $\rho_{\alpha, t}$ is the pointwise limit of the sequence $(\rho_{\alpha, t} ^{(j)}) _{j \in \mathbb N}$ corresponding to an exhaustion of $M$ with regular domains, we could have first obtained the function $\operatorname{Strat}_t ^{(j)} (\alpha)$ on $\mathcal C_t (\overline {U_j})$ for every $j \in \mathbb N$, shown that these satisfy a natural compatibility relation, and then defined $\operatorname{Strat}_t (\alpha)$ as the pointwise limit of these functions. The drawback of such an approach is that proving the convergence $S_{P,t,k} (\alpha) \to \operatorname{Strat} _t (\alpha)$ becomes more difficult, and proving the uniformity with respect to $t$ of this convergence becomes even more difficult. Therefore, the approach that we have finally opted for in this article has been motivated by the desire for simplicity of the proofs involved.

\section{A general concept of stochastic integral}

In order to unravel a general concept of stochastic integral, let us return to the approximations $S_{P,t,k} (\alpha)$ constructed above and define the related approximations
\begin{align}
\label{approximants} A_{P,t,k} (\alpha) (c) & = S_{P,t,k} (\alpha) (c) - \frac t {2^k} \sum _{j=0} ^{2^k - 1} (\mathrm d ^* \alpha) \left( c \left( \frac {jt} {2^k} \right) \right) \int _{[0,1]} (2 \tau - 1) \, \mathrm d P (\tau) = \\
\nonumber & = \sum _{j=0} ^{2^k - 1} I_P (\alpha) \left( c \left( \frac {jt} {2^k} \right), c \left( \frac {(j+1)t} {2^k} \right) \right)
\end{align}
for every curve $c \in \mathcal C_t$ (that is, we just drop the term containing $\mathrm d ^* \alpha$). We shall now study the behaviour of these approximations on continuously differentiable curves, this "classical" behaviour going to guide us towards the understanding of its "stochastic" counterpart.

\subsection{An approximation of the line integral on differentiable curves}

In order to shorten the formulae, from now on we shall use the notation $p_{k,j} = c \left( \dfrac {jt} {2^k} \right)$ for every continuous $c : [0,t] \to M$, every $k \in \mathbb N$ and every $0 \le j \le 2^k$.

The next lemma has a purely technical and auxiliary character, its role being to allow us to control the geometry in the neighbourhood of some continuously differentiable curve $c$, which in turn will help us obtain majorizations indispensable to the proof of theorem \ref{classical limit}.

\begin{lemma} \label{preliminaries}
If $c : [0,t] \to M$ is a continuously differentiable curve, then there exists $k_c \in \mathbb N$ and a compact subset $K \subseteq M$ such that:
\begin{enumerate}[wide]
\item For every $k \ge k_c$ and every $0 \le j \le 2^k - 1$, the points $p_{k,j}$ and $p_{k, j+1}$ may be joined by a unique minimizing geodesic.

\item The image of the curve $c$ is contained in $K$.

\item If $\gamma_{k,j} : [0,1] \to M$ is the unique minimizing geodesic joining $p_{k,j}$ to $p_{k,j+1}$, that is $\gamma_{k,j} (s) = \exp_{p_{k,j}} (s \exp_{p_{k,j}} ^{-1} (p_{k,j+1}))$, and if $\gamma_k : [0,t] \to M$ is the polygonal line formed by joining the segments $\gamma_{k,j}$ end to end, that is $\gamma_k (s) = \gamma_{k,j} ( \frac {2^k} t s - j)$ for $s \in \left[ \dfrac {jt} {2^k}, \dfrac {(j+1)t} {2^k} \right]$, then the image of the polygonal line $\gamma_k$ is contained in $K$ for every $k \ge k_c + 1$.
\end{enumerate}
\end{lemma}
\begin{proof}
\begin{enumerate}[wide]
\item Since $c([0,t])$ is a compact subset, $r_c = \frac 1 2 \min(1, \min _{s \in [0,t]} \operatorname{injrad} (c(s)))$ exists and is non-zero. (The truncation at $1$ ensures that $r_c$ is finite, whereas the factor $\frac 1 2$ could be replaced by any other number from $(0,1)$, its sole purpose being to guarantee the compactness of the closed ball $\overline B (c(s), r_c) = \{x \in M \mid d(x, c(s)) \le r_c\}$, this being diffeomorphic under $\exp_{c(s)} ^{-1}$ to a closed ball of $T_{c(s)} M$ for every $s \in [0,t]$.)

Let us consider two consecutive points $p_{k,j}$ and $p_{k,j+1}$ on $c$. We are looking for a sufficient condition allowing for them to be joined by a unique minimizing geodesic. We impose, therefore, that $d(p_{k,j}, p_{k,j+1}) \le r_c$, whence $\frac t {2^k} \| \dot c \| _\infty \le r_c$, so if $k_c = \left\lceil \log_2 \frac {t \, \| \dot c \| _\infty} {r_c} \right\rceil$ then $p_{k,j}$ and $p_{k,j+1}$ may be joined by a unique minimizing geodesic for every $k \ge k_c$ and $0 \le j \le 2^k - 1$.

\item Let $K = \bigcup _{j = 0} ^{2^{k_c + 1} - 1} \overline B (p_{k_c + 1,j}, r_c)$. For every $k \in \mathbb N$ and $0 \le j \le 2^k - 1$, let $c_{k,j} : [0,1] \to M$ be the restriction of the curve $c$ to the interval $\left[ \dfrac {jt} {2^k}, \dfrac {(j+1)t} {2^k} \right]$, that is $c_{k,j} (s) = c \left( \dfrac {(j+s)t} {2^k} \right)$. With these notations, we shall show that the image of the curve $c_{k_c + 1,j}$ is contained in $\overline B (p_{k_c + 1,j}, r_c)$ for every $0 \le j \le 2^{k_c + 1}$. Indeed, if $0 \le s \le 1$, then
\[ d (p_{k_c + 1,j}, c_{k_c,j} (s)) = d \left( c \left( \frac {jt} {2^{k_c + 1}} \right), c \left( \frac {(j+s)t} {2^{k_c + 1}} \right) \right) \le \frac {st} {2^{k_c + 1}} \| \dot c \| _\infty \le \frac t {2^{k_c + 1}} \| \dot c \| _\infty \le \frac {r_c} 2 < r_c \ . \]
Since $c$ is the union of the curves $c_{k_c, j}$ when $0 \le j \le 2^{k_c}-1$, it follows that the image of the curve $c$ is contained in $K$.

\item If $k \ge k_c + 1$, and $0 \le j_0 \le 2^{k_c + 1} - 1$ is such that $\dfrac {j_0} {2^{k_c + 1}} \le \dfrac j {2^k} \le \dfrac {j_0 + 1} {2^{k_c + 1}}$, then for every $s \in [0,t]$ we have
\begin{align*}
d(p_{k_c + 1, j_0}, \gamma_{k,j}(s)) & \le d (p_{k_c + 1,j_0}, p_{k,j}) + d (p_{k,j}, \gamma_{k,j} (s)) \le \\
& \le \| \dot c \| _\infty \left( \frac {jt} {2^k} - \frac {j_0 t} {2^{k_c + 1}} \right) + \| \dot c \| _\infty \, \frac s {2^k} \le \| \dot c \| _\infty \, \frac t {2^{k_c + 1}} + \| \dot c \| _\infty \, \frac t {2^{k_c + 1}} \le r_c \ ,
\end{align*}
so the image of $\gamma_{k,j}$ is contained in $\overline B(p_{k_c + 1,j}, r_c)$, whence the image of $\gamma_k$ is contained in $K$.
\end{enumerate}
\end{proof}

\begin{theorem} \label{classical limit}
If $c : [0,t] \to M$ is a twice continuously differentiable curve, then
\[ \int _c \alpha = \lim _{k \to \infty} A_{P,t,k} (\alpha) (c) \ . \]
\end{theorem}

\begin{proof}
We shall use the notations introduced in lemma \ref{preliminaries} and its proof. Given that we are interested in a limit when $k \to \infty$, we shall work under the hypothesis that $k \ge k_c + 1$, so that we may benefit from all the conclusions in that lemma. Also, if $f$ is some continuous function defined on $K$, we shall denote the maximum of its absolute value by $\| f \| _{\infty, K}$. 

We are looking, then, to estimate the difference
\[ \left| A_{P,t,k} (\alpha) (c) - \int _c \alpha \right| = \left| \sum _{j=0} ^{2^k - 1} \int _{[0,1]} \alpha _{\gamma_{k,j} (\tau)} ( \dot \gamma_{k,j} (\tau) ) \, \mathrm d P (\tau) - \sum _{j=0} ^{2^k - 1} \int _{\frac {jt} {2^k}} ^{\frac {(j+1)t} {2^k}} \alpha _{c(\theta)} ( \dot c(\theta) ) \, \mathrm d \theta \right| \ . \]

Since $\gamma_{k,j} (0) = p_{k,j}$ and $\dot \gamma_{k,j} (0) = \exp _{p_{k,j}} ^{-1} (p_{k,j+1})$, and since $\nabla _{\dot \gamma_{k,j}} \dot \gamma_{k,j} = 0$ (because $\gamma_{k,j}$ is a geodesic), the first integral is
\begin{align*}
\int _{[0,1]} & \alpha _{\gamma_{k,j} (\tau)} ( \dot \gamma_{k,j} (\tau) ) \, \mathrm d P (\tau) = \\
& = \int _{[0,1]} \left( \alpha _{p_{k,j}} ( \exp _{p_{k,j}} ^{-1} (p_{k,j+1}) ) + \int _{[0,\tau]} (\nabla \alpha) _{\gamma_{k,j} (\sigma)} (\dot \gamma_{k,j} (\sigma), \dot \gamma_{k,j} (\sigma)) \, \mathrm d \sigma \right) \mathrm d P (\tau) \ ,
\end{align*}
where $\nabla$ is the Levi-Civita connection acting on $1$-forms. The first term in the integral may be expanded using Taylor's formula up to order $2$,
\begin{align*}
\alpha_{p_{k,j}} & \left( \exp _{p_{k,j}} ^{-1} (p_{k,j+1}) \right) = \alpha_{p_{k,j}} \left( \exp _{c \left( \frac {jt} {2^k} \right)} ^{-1} \left( c\left( \frac {(j+1)t} {2^k} \right) \right) \right) = \\
& = \alpha_{p_{k,j}} \left( \dot c \left( \frac {jt} {2^k} \right) \right) \frac t {2^k} + \int _{\frac {jt} {2^k}} ^{\frac {(j+1)t} {2^k}} \left( \frac {(j+1)t} {2^k} - \theta \right) \, \alpha_{p_{k,j}} \left( \frac {\mathrm d ^2} {\mathrm d \theta^2} \exp _{c \left( \frac {jt} {2^k} \right)} ^{-1} (c(\theta)) \right) \mathrm d \theta \ .
\end{align*}

Similarly, using Taylor's formula up to order $1$, the second term in the integral is
\begin{align*}
& \int _{\frac {jt} {2^k}} ^{\frac {(j+1)t} {2^k}} \alpha _{c(\theta)} ( \dot c(\theta) ) \, \mathrm d \theta = \\
& = \alpha_{p_{k,j}} \left( \dot c \left( \frac {jt} {2^k} \right) \right) \frac t {2^k} + \int _{\frac {jt} {2^k}} ^{\frac {(j+1)t} {2^k}} \int _{\frac {jt} {2^k}} ^\theta \alpha _{c(\sigma)} ( (\nabla \alpha) _{c(\sigma)} (\dot c (\sigma), \dot c (\sigma)) + \nabla _{\dot c (\sigma)} \dot c (\sigma)) \, \mathrm d \sigma \, \mathrm d \theta \ .
\end{align*}

Combining the results obtained so far, and taking into account the fact that $\int _{[0,1]} \mathrm d P(\tau) = 1$, we obtain that
\begin{align*}
\left| A_{P,t,k} (\alpha) (c) - \int _c \alpha \right| & \le \sum _{j=0} ^{2^k - 1} \left| \int _{\frac {jt} {2^k}} ^{\frac {(j+1)t} {2^k}} \left( \frac {(j+1)t} {2^k} - \theta \right) \, \alpha_{p_{k,j}} \left( \frac {\mathrm d ^2} {\mathrm d \theta^2} \exp _{c \left( \frac {jt} {2^k} \right)} ^{-1} (c(\theta)) \right) \mathrm d \theta \right| + \\
& + \sum _{j=0} ^{2^k - 1} \left| \int _{[0,1]} \int _{[0,\tau]} (\nabla \alpha) _{\gamma_{k,j} (\sigma)} (\dot \gamma_{k,j} (\sigma), \dot \gamma_{k,j} (\sigma)) \, \mathrm d \sigma \, \mathrm d P (\tau) \right| + \\
& + \sum _{j=0} ^{2^k - 1} \left| \int _{\frac {jt} {2^k}} ^{\frac {(j+1)t} {2^k}} \int _{\frac {jt} {2^k}} ^\tau \alpha _{c(\sigma)} ( (\nabla \alpha) _{c(\sigma)} (\dot c (\sigma), \dot c (\sigma)) + \nabla _{\dot c (\sigma)} \dot c (\sigma)) \, \mathrm d \sigma \, \mathrm d \tau \right| .
\end{align*}
We shall majorize each of these three sums, one by one.

In the integral in the first sum we notice that $\theta - \frac {jt} {2^k} \le \frac t {2^k} \le \frac t {2^{k_c + 1}}$ (remember that we have chosen $k_c$ such that any two points $c(s)$ and $c(s')$ may be joined by a unique minimizing geodesic if $|s-s'| \le \frac t {2^{k_c}}$), so
\begin{gather*}
\sup \left\{ \left\| \frac {\mathrm d ^2} {\mathrm d \theta^2} \exp _{c \left( \frac {jt} {2^k} \right)} ^{-1} (c(\theta)) \right\| _{T_{p_{k,j}} M} \mid k \ge k_c, \, j \in \{0, \dots, 2^k-1\}, \, \theta \in \left[ \frac {jt} {2^k}, \frac {(j+1)t} {2^k} \right] \right\} \le \\
\le \sup \left\{ \left\| \frac {\mathrm d ^2} {\mathrm d \theta^2} \exp _{c (s)} ^{-1} (c(s+\theta)) \right\| _{T_{p_{k,j}} M} \mid s \in [0,t], \, \theta \in \left[0, \frac 1 {2^{k_c}} \right] \right\} \ne \infty
\end{gather*}
because the Riemannian exponential map is smooth in both arguments, $c$ is twice continuously differentiable and $[0,t] \times \left[0, \frac 1 {2^{k_c}} \right]$ is compact. We conclude that there exists $C_1 > 0$ such that the first sum is at most
\[ C_1 \sup _{s \in [0,t]} \| \alpha _{c(s)} \| _{T^* _{c(s)} M} \sum _{j=0} ^{2^k - 1} \int _{\frac {jt} {2^k}} ^{\frac {(j+1)t} {2^k}} \left( \frac {(j+1)t} {2^k} - \theta \right) \mathrm d \theta = \frac 1 2 C_1 \| \alpha \| _{\infty, K} \, \frac {t^2} {2^k} \]
when $k \ge k_c$.

For the second sum notice that the curves $\gamma_{k,j}$ have constant length tangent vectors (being geodesics), and remember that the Riemannian exponential map is a radial isometry, which leads to
\[ \| \dot \gamma _{k,j} (\sigma) \| _{T _{\gamma_{k,j} (\sigma)} M} = \| \exp _{p_{k,j}} ^{-1} (p_{k,j+1})\| = d(p_{k,j}, p_{k,j+1}) \le \int _{\frac {jt} {2^k}} ^{\frac {(j+1)t} {2^k}} \| \dot c (s) \| _{T_{c(s)} M} \le \| \dot c \| _{\infty, [0,t]} \, \frac t {2^k} \ ,\]
whence that sum may be majorized by
\[ \frac 1 2 \, \| \nabla \alpha \| _{\infty, K} \, \| \dot c \| _{\infty, [0,t]} ^2 \, \frac {t^2} {2^k} \]
when $k \ge k_c + 1$.

Finally, the third sum may be majorized by
\[ \frac 1 2 \| \alpha \| _{\infty, K} \, (\| \dot c \| _{\infty, [0,t]} ^2 + \| \nabla _{\dot c} \dot c \| _{\infty, [0,t]}) \, \frac {t^2} {2^k} \ .\]

We conclude that there exists $C>0$ (depending on $c$, $\alpha$ and the compact $K$) such that
\[ \left| A_{P,t,k} (\alpha) (c) - \int _c \alpha \right| \le C \, t^2 \, \frac 1 {2^k} \]
for $k \ge k_c + 1$, whence the conclusion.
\end{proof}

\subsection{A geometrical definition and a classification of stochastic integrals}

We shall draw inspiration from the resemblance between theorem \ref{classical limit} and theorem \ref{approximation of the Stratonovich integral} in order to exhibit a general concept of stochastic integral. Let $\operatorname{Prob} ([0,1])$ be the space of regular Borel probability measures on the interval $[0,1]$.

\begin{definition}
We shall say that $\operatorname{Int}_t : \Omega ^1 (M) \to L^0 (\mathcal C_t)$ is a stochastic integral if and only if there exists $P \in \operatorname{Prob} ([0,1])$ such that $\operatorname{Int}_t (\alpha)$ be the limit in measure of the sequence of approximations $A_{P,t,k} (\alpha)$ for every $\alpha \in \Omega^1 (M)$. When this condition is met, we shall denote this stochastic integral by $\operatorname{Int}_{P,t}$, in order to emphasize its dependence on $P$.
\end{definition}

Although the convergence in measure obtained in theorem \ref{approximation of the Stratonovich integral} was uniform with respect to $t$ in bounded subsets of $(0, \infty)$, we have not included this property in the above definition because it was not clear, upon writing this text, whether this uniformity is an essential ingredient of the concept or a merely accidental one without major consequences.

\begin{remark}
Given that convergence in measure (as in the proposed definition) is weaker than pointwise convergence, let us emphasize that $\operatorname{Int}_t (\alpha) (\cdot)$ must be understood not as a function defined for every curve from $\mathcal C_t$, but rather as an element from $L^0 (\mathcal C_t)$. This is the major difference from the usual line integral, which is defined for every piecewise-differentiable curve.
\end{remark}

Let $P \in \operatorname{Prob} ([0,1])$. We would like to discover whether there exists any connection between the freshly defined $\operatorname{Int}_{P,t} (\alpha)$ and the function $\operatorname{Strat} _t (\alpha)$ obtained in \ref{the Stratonovich integral}. Let us notice that
\[ \lim _{k \to \infty} \frac t {2^k} \sum _{j=0} ^{2^k - 1} (\mathrm d ^* \alpha) \left( c \left( \frac {jt} {2^k} \right) \right) = \int _0 ^t (\mathrm d^* \alpha) (c(s)) \, \mathrm d s \]
for every $c \in \mathcal C_t$, as the limit of the Riemann sums associated to the continuous function $(\mathrm d^* \alpha) \circ c$, the equidistant partition of $[0,t]$ into $2^k$ subintervals, and the intermediate points $\big( c (\frac {jt} {2^k}) \big) _{0 \le j \le 2^k - 1}$. Even more, then, is the above convergence also valid in measure. If we pass to the limit in measure in formula \eqref{approximants} used to define the approximations $A_{P,t,k}$, we get
\[ \operatorname{Int}_{P,t} (\alpha) (c) = \operatorname{Strat} _t (\alpha) (c) - \int _{[0,1]} (2 \tau - 1) \, \mathrm d P (\tau) \int _0 ^t (\mathrm d^* \alpha) (c(s)) \, \mathrm d s \ , \]
which shows that although the probability $P$ may be extremely complicated, the corresponding stochastic integral $\operatorname{Int}_{P,t}$ remembers only its first-order moment, discarding any other information associated to $P$; furthermore, any two probabilities from $\operatorname{Prob} ([0,1])$ with the same first order moment give rise to the same stochastic integral. We also conclude that, since the function $\operatorname{Strat} _t$ has already been constructed, $\operatorname{Int}_{P,t}$ exists for every $P \in \operatorname{Prob} ([0,1])$. Since the function $2\tau - 1$ has the minimum $-1$ and the maximum $+1$ on $[0,1]$, and since $P$ is a probability, it follows that $\int _{[0,1]} (2\tau - 1) \, \mathrm d P (\tau) \in [-1, 1]$, and that every stochastic integral $\operatorname{Int}_t$ on $\mathcal C_t$ is of the form
\[ \operatorname{Int}_t (\alpha) = \operatorname{Strat} _t (\alpha) + \theta \int _0 ^t (\mathrm d ^* \alpha) (c(s)) \, \mathrm d s \]
with $\theta \in [-1, 1]$.

Furthermore, if $P,Q \in \operatorname{Prob} ([0,1])$, then
\[ \operatorname{Int}_{P,t} (\alpha) (c) = \operatorname{Int}_{Q,t} (\alpha) (c) - 2 \int _{[0,1]} \tau \, \mathrm d (P-Q) (\tau) \int _0 ^t (\mathrm d^* \alpha) (c(s)) \, \mathrm d s \ ,\]
so that any two stochastic integrals differ by a multiple of the integral of $\mathrm d^* \alpha$.

This is a good moment to see several concrete examples of such stochastic integrals as defined in this work, and to compare our results to the ones already obtained in the stochastic literature.

\begin{itemize}[wide]
\item If $P = \delta_0$ (the Dirac measure concentrated at $0$), then
\[ \operatorname{Int}_{\delta_0, t} (\alpha) (c) = \operatorname{Strat} _t (\alpha) (c) + \int _0 ^t (\mathrm d^* \alpha) (c(s)) \, \mathrm d s \ . \]
By comparing our approximations $A_{\delta_0, k, t} (\alpha) $ of $I_{\delta_0, t} (\alpha)$ to the ones in theorem 7.37 on page 110 of \cite{Emery89} (or to the ones in theorem A from \cite{Darling84}, which is nevertheless stated under more restrictive hypotheses than here), we recognize immediately that $\operatorname{Int}_{\delta_0, t} (\alpha)$ is the Itô integral of $\alpha$, therefore from now on we shall denote it by $\operatorname{Ito}_t (\alpha)$.

\item If $P = \operatorname{Leb}_{[0,1]}$ (the Lebesgue measure on $[0,1]$), or $P = \delta _{\frac 1 2}$ (the Dirac measure concentrated at $\frac 1 2$), or $P = \frac 1 2 \delta _1$, or $P = \frac 1 2 (\delta_0 + \delta_1)$, then the corresponding stochastic integral is
\[ \operatorname{Int}_{\operatorname{Leb}_{[0,1]}, t} (\alpha) = \operatorname{Strat} _t (\alpha) \ . \]
By comparing the approximations $A_{\operatorname{Leb}_{[0,1]}, k, t} (\alpha)$ of $\operatorname{Int}_{\operatorname{Leb}_{[0,1]}, t} (\alpha)$ to those in theorem 7.14 on page 96 of \cite{Emery89}, we readily recognize that $\operatorname{Int}_{\operatorname{Leb}_{[0,1]}, t} (\alpha)$ is the Stratonovich integral of $\alpha$. (The reader is invited to compare these results to the ones in section 6 of \cite{Norris92}, too.)

\item In general, if $M_1 (P) = \int _{[0,1]} \tau \, \mathrm d P (\tau)$, then the stochastic integral $\operatorname{Int}_{P,t} (\alpha)$ corresponding to $P \in \operatorname{Prob}([0,1])$ coincides with the one produced by the probabilities $\delta _{M_1 (P)}$ (the Dirac measure concentrated at $M_1 (P) \in [0,1]$) and $(1 - M_1 (P)) \delta_0 + M_1 (P) \delta_1$, all these probabilities having $M_1 (P)$ as first order moment. Nevertheless, although in principle we could study the stochastic integrals defined in this work using only these very simple combinations of Dirac measures, some results are much easier to prove using more complicated probabilities with the same first order moment. In particular, in the study of the Stratonovich integral it is usually more convenient to use the Lebesgue measure on $[0,1]$.

\item With this new insight into the problem, it is worth reviewing now the construction of the operator that was given by
\[ (R_{\alpha, t} ^{(U)} f) (x) = \int _U h_\alpha ^{(U)} (t,x,y) \, \chi(x,y) \, \mathrm e ^{-\mathrm i I_P (\alpha) (x,y) - \mathrm i t (\mathrm d^* \alpha) (x) \int _{[0,1]} (2 \tau - 1) \, \mathrm d P (\tau)} f(y) \, \mathrm d y \ . \]
The term containing $\mathrm d ^* \alpha$ seems out of place and artificial; indeed, it was added only as a correction term, in the absence of which theorem \ref{the technical theorem} would not have held anymore. In turn, the strength of theorem \ref{the technical theorem} lies in the fact that the left-hand side depends on $P$, while the right-hand side does not. This is the crucial ingredient that allows us to link any two stochastic integrals (as defined above) by a simple formula, showing that in fact there exists essentially a single stochastic integral (any of them would do). If, on the other hand, we were to concentrate only on the study of the Stratonovich integral (which, as it had become apparent in this work, plays a special role among the stochastic integrals), then we could choose $P = \operatorname{Leb}_{[0,1]}$, in which case the operator $R_{\alpha, t} ^{(U)}$ would get the much simpler form
\[ (R_{\alpha, t} ^{(U)} f) (x) = \int _U h_\alpha ^{(U)} (t,x,y) \, \chi(x,y) \, \mathrm e ^{-\mathrm i I_{\operatorname{Leb}_{[0,1]}} (\alpha) (x,y)} f(y) \, \mathrm d y \ . \]
Disregarding the cut-off function which is there only for technical reasons, the above formula says to take the vector $f(y)$ from the fiber $\{y\} \times \mathbb C$ of the trivial bundle $M \times \mathbb C$, parallel-transport it with respect to the connection $\mathrm d - \mathrm i \alpha$ along the geodesic $\gamma_{y,x}$ from $y$ to $x$, multiply by the heat kernel and integrate. It is worth noting that this purely-geometric recipe could be used verbatim in any other Hermitian bundle endowed with a Hermitian connection over $M$, which is expected to lead to the obtention of the stochastic parallel transport. This line of thought will not be followed anymore in this work, though, but will form the subject of a future one.
\end{itemize}

Let us quickly discuss the orientation of the tangent vectors in the approximations of the Itô integral: why into the future and not into the past? If we understand stochastic integration in the way presented here, the orientation of the tangent vectors $\exp_x ^{-1} y$ (or $y-x$ in $\mathbb R^n$) is baked into the expression $\int _{[0,1]} \alpha _{\gamma_{x,y} (\tau)} (\dot \gamma_{x,y} (\tau)) \, \mathrm d P (\tau)$ (with $P = \delta_0$), that is into the expression $\alpha _{c(s)} (\dot c (s))$ that appears in the usual line integral. We therefore understand that the orientation of the tangent vectors is not something that has to be \textit{chosen}, but rather an automatic, inevitable, consequence of the analogy between stochastic integrals and line integrals that has been guiding our intuition from the very beginning.

\begin{remark}
The previous examples show that the Stratonovich and Itô integrals of $\alpha$ are equal if and only if $\mathrm d^* \alpha = 0$. This is worth comparing to lemma 8.24 on page 120 of \cite{Emery89}, where only a necessary but not sufficient condition (difficult to verify in concrete applications) is given that guarantees this equality. More specifically, Émery first introduces the concept of stochastic parallel transport in the bundles $TM$ and $T^*M$, starting from which he constructs certain martingales depending on $\alpha$; if these martingales are of finite variation, then the Stratonovich and Itô integrals of $\alpha$ are equal.
\end{remark}

\section{Basic properties of stochastic integrals}

We are now going to see that the theoretical construction that we have performed is useful and pays off. The following two results are fundamental in stochastic analysis, but in our functional-analytic approach of the subject they are just elementary and natural consequences of everything done so far.

\begin{proposition}
The Stratonovich integral $\operatorname{Strat} _t$ has the property that $\operatorname{Strat} _t (\mathrm d f) (c) = f(c(t)) - f(x_0)$ for every real-valued continuously-differentiable function $f$ and for almost every $c \in \mathcal C_t$.
\end{proposition}
\begin{proof}
Using the approximation of the Stratonovich integral based upon the Lebesgue measure,
\[ \operatorname{Strat} _t (\mathrm d f) (c) = \lim_{k \to \infty} \sum _{j=0} ^{2^k - 1} I_{\operatorname{Leb}_{[0,1]}, t} (\mathrm d f) \left( c \left( \frac {jt} {2^k} \right), c \left( \frac {(j+1)t} {2^k} \right) \right) \]
in measure, whence there exists a sub-sequence $(k_l) _{l \ge 0} \subseteq \mathbb N$ such that
\[ \operatorname{Strat} _t (\mathrm d f) (c) = \lim_{l \to \infty} \sum _{j=0} ^{2^{k_l} - 1} I_{\operatorname{Leb}_{[0,1]}, t} (\mathrm d f) \left( c \left( \frac {jt} {2^{k_l}} \right), c \left( \frac {(j+1)t} {2^{k_l}} \right) \right) \]
for almost every $c \in \mathcal C_t$. Denote by $\mathcal C'$ the co-null subset of $\mathcal C_t$ on which this equality is true.

Let $c \in \mathcal C'$. Since $[0,t]$ is compact, $c$ will be uniformly continuous. Let $\omega$ be a non-decreasing modulus of continuity for $c$, that is $d(c(s), c(s')) \le \omega(|s-s'|)$ for every $s,s' \in [0,t]$, with $\omega$ non-decreasing and $\omega(0) = 0$. (We could use here that the Wiener measure is concentrated on the Hölder-continuous curves, which have an explicit modulus of continuity, but we shall pretend that we do not know this fact.) Let us show that, from a certain $k_c \in \mathbb N$ onward, any two consecutive points $p_{k_l, j}$ and $p_{k_l, j+1}$ on $c$ may be joined by a unique minimizing geodesic. Indeed, since $c([0,t])$ is compact, there exists $r_c = \frac 1 2 \min (1, \min _{s \in [0,t]} \operatorname{injrad} (c(s)) ) \in (0, \infty)$. Therefore, if we require that
\[ d(p_{k_l, j}, p_{k_l, j+1}) \le \omega \left( \frac t {2^{k_l}} \right) < r_c \ , \]
given that $\omega$ is non-decreasing and $\omega(0) = 0$, it follows that there exists $k_c \in \mathbb N$ such that $p_{k_l, j}$ and $p_{k_l, j+1}$ may be joined by a unique minimizing geodesic for every $k_l \ge k_c$ and $0 \le j \le 2^{k_l} - 1$. Since $k_l \ge k_c$ when $l \to \infty$, we have that
\begin{align*}
\operatorname{Strat} _t (\mathrm d f) (c) & = \lim_{l \to \infty} \sum _{j=0} ^{2^{k_l} - 1} I_{\operatorname{Leb}_{[0,1]}, t} (\mathrm d f) \left( c \left( \frac {jt} {2^{k_l}} \right), c \left( \frac {(j+1)t} {2^{k_l}} \right) \right) = \\
& = \lim_{l \to \infty} \sum _{j=0} ^{2^{k_l} - 1} \int _{[0,1]} \mathrm d f _{\gamma_{p_{k_l,j},p_{k_l,j+1}} (\tau)} (\dot \gamma_{p_{k_l,j},p_{k_l,j+1}} (\tau)) \, \mathrm d \tau = \\
& = \lim_{l \to \infty} \sum _{j=0} ^{2^{k_l} - 1} \int _0 ^1 \frac {\mathrm d} {\mathrm d \tau} (f \circ \gamma_{p_{k_l,j},p_{k_l,j+1}}) (\tau) \, \mathrm d \tau = \\
& = \lim_{l \to \infty} \sum _{j=0} ^{2^{k_l} - 1} [f (p_{k_l, j+1}) - f (p_{k_l, j+1})] = \\
& = f(c(t)) - f(x_0) \ ,
\end{align*}
for every $c \in \mathcal C'$, and since $\mathcal C' \subseteq \mathcal C_t$ is co-null, the conclusion is proved.
\end{proof}

\begin{corollary}[Itô's lemma]
If $f : M \to \mathbb R$ is real, twice continuously-differentiable, and $\Delta$ is the Laplace-Beltrami operator on $M$, then
\[ f(c(t)) = f(x_0) + \operatorname{Ito}_t (\mathrm d f) (c) + \int _0 ^t (\Delta f) (c(s)) \, \mathrm d s \]
for almost every $c \in \mathcal C_t$.
\end{corollary}
\begin{proof}
Using the previous theorem, the proof is short and elementary:
\begin{align*}
f(c(t)) & = f(x_0) + \operatorname{Strat} _t (\mathrm d f) (c) = f(x_0) + \operatorname{Ito}_t (\mathrm d f) (c) - \int _0 ^t \mathrm d ^* (\mathrm d f) (c(s)) \, \mathrm d s = \\
& = f(x_0) + \operatorname{Ito}_t (\mathrm d f) (c) + \int _0 ^t (\Delta f) (c(s)) \, \mathrm d s \ ,
\end{align*}
where we have used the basic Hodge-theoretical formula $\Delta = - \mathrm d ^* \mathrm d$.
\end{proof}

\begin{remark}
The reason why some authors obtain a $\frac 1 2$ factor in front of the Laplacian is their use of the convention in which the heat operator is defined as $\partial_t - \frac 1 2 \Delta$, whereas in the present work the heat operator is $\partial_t - \Delta$.
\end{remark}

As always when one studies objects that depend on certain parameters, it is useful to study how regular this dependency is. In particular, it is interesting to study the dependence of the stochastic integral $\operatorname{Int}_{P,t} (\alpha)$ on the parameter $t \in (0,T]$, where $T>0$ is arbitrary. Since the stochastic integral $\operatorname{Int}_{P,t} (\alpha)$ lives in the space $L^0 (\mathcal C _t)$ for each $t \in (0, T]$, and since all these spaces are unrelated to each other, we shall have to embed all of them in the bigger space $L^0 (\mathcal C _T)$. In order to do this, let us remember that the natural topology in $L^0 (\mathcal C _t)$ is that of convergence in the Wiener measure $w_t$. If $\operatorname{res}_{[0,t]} : \mathcal C _T \to \mathcal C _t$ is the restriction $\operatorname{res}_{[0,t]} (c) = c | _{[0,t]}$, then clearly $w_t = (\operatorname{res}_{[0,t]} ) _* w_T$. This topology is metrizable by any distance of the form 
\begin{gather*}
d_t (f,g) = \int _{\mathcal C _t} \varphi (|f-g|) \, \mathrm d w_t = \\
= \int _{\mathcal C _T} \varphi (|f \circ \operatorname{res}_{[0,t]} - g \circ \operatorname{res}_{[0,t]}|) \, \mathrm d w_T = d_T (f \circ \operatorname{res}_{[0,t]}, g \circ \operatorname{res}_{[0,t]})
\end{gather*}
with $\varphi : [0, \infty) \to [0, \infty)$ continuous, bounded, concave, non-decreasing, with $\varphi(0) = 0$ and $\varphi > 0$ on $(0, \infty)$. Any such distance is called a "Lévy distance".

For the line integral, if $c : [0,T] \to M$ is continuously-differentiable, then
\[ \left| \int _c \alpha - \int _{\operatorname{res}_{[0,t]} (c)} \alpha \right| = \left| \int _{\operatorname{res}_{[t,T]} (c)} \alpha \right| \le \sup _{s \in [0,T]} | \alpha _{c(s)} (\dot c (s)) | \, (T-t) \ . \]
In particular, the map $[0,T] \ni t \mapsto \int _{\operatorname{res}_{[0,t]} (c)} \alpha \in \mathbb R$ is continuous. The following theorem offers a weaker analogue of this fact in the context of stochastic integration.
 
\begin{theorem}
For every $\alpha \in \Omega^1 (M)$, the map $(0,T] \ni t \mapsto \operatorname{Int}_{P,t} (\alpha) \circ \operatorname{res}_{[0,t]} \in L^0 (\mathcal C _T)$ is continuous.
\end{theorem}
\begin{proof}
We want to show that if $t_m \to t$ in $(0,T]$, then $\operatorname{Int}_{P,t_m} (\alpha) \circ \operatorname{res}_{[0,t_m]} \to \operatorname{Int}_{P,t} (\alpha) \circ \operatorname{res}_{[0,t]}$ in $L^0 (\mathcal C_T, w_T)$. Fix a function $\varphi$ as described above and let $d_t$ be the corresponding Lévy distance that it generates on the space $L^0 (\mathcal C_T)$, for all $t \in (0,T]$ (the same $\varphi$ is used for all $t$); using the triangle inequality gives us that
\begin{gather*}
d_T (\operatorname{Int}_{P,t_m} (\alpha) \circ \operatorname{res}_{[0,t_m]}, \operatorname{Int}_{P,t} (\alpha) \circ \operatorname{res}_{[0,t]}) \le \\
\le  d_T (\operatorname{Int}_{P,t_m} (\alpha) \circ \operatorname{res}_{[0,t_m]}, A_{P,t_m,k} (\alpha) \circ \operatorname{res}_{[0,t_m]}) +  \\
+ d_T (A_{P,t_m,k} (\alpha) \circ \operatorname{res}_{[0,t_m]}, A_{P,t,k} (\alpha) \circ \operatorname{res}_{[0,t]}) + \\
+ d_T (A_{P,t,k} (\alpha) \circ \operatorname{res}_{[0,t]}, \operatorname{Int}_{P,t} (\alpha) \circ \operatorname{res}_{[0,t]}) = \\
= d_{t_m} (\operatorname{Int}_{P,t_m} (\alpha), A_{P,t_m,k} (\alpha)) + d_T (A_{P,t_m,k} (\alpha) \circ \operatorname{res}_{[0,t_m]}, A_{P,t,k} (\alpha) \circ \operatorname{res}_{[0,t]}) + \\
+ d_t (A_{P,t,k} (\alpha), \operatorname{Int}_{P,t} (\alpha))
\end{gather*}

Since we have $A_{P,t,k} (\alpha) \to \operatorname{Int}_{P,t} (\alpha)$ in the measure $w_t$ by our very definition of stochastic integrals (and the same for $t_m$ instead of $t$), we conclude that the first and third terms converge to $0$.

For the middle term, we shall show that $A_{P,t_m,k} (\alpha) (c | _{[0,t_m]}) \to A_{P,t,k} (\alpha) (c | _{[0,t]})$ when $t_m \to t$, for every $c \in \mathcal C _T$, from which the weaker convergence in the measure $w_T$ will follow. We want to show, thus, that
\[ \sum _{j=0} ^{2^k-1} I_P (\alpha) \left( c \left( \frac {j t_m} {2^k} \right), c \left( \frac {(j+1) t_m} {2^k} \right) \right) \to \sum _{j=0} ^{2^k-1} I_P (\alpha) \left( c \left( \frac {j t} {2^k} \right), c \left( \frac {(j+1) t} {2^k} \right) \right) \]
but, for any sufficiently large $k$, this is an elementary consequence of the continuity of $c$ and of $I_P (\alpha)$.

Since each of the three terms in the right hand side tends to $0$, the left hand side will also tend to $0$.
\end{proof}

\section{An application: the Feynman-Kac-Itô formula}

In the following theorem we shall see that the Stratonovich integral emerges absolutely naturally when we try to deduce the analogue of the Feynman-Kac formula in the presence of a magnetic field represented by the $1$-form $\alpha$. Theorem 15.3 on page 162 of \cite{Simon79} presents an alternative view on the same problem, but only for the particular case $M = \mathbb R^n$. For every $x \in M$ we shall consider the space
\[ \mathcal C _{t,x} = \{ c : [0,t] \to M \mid c \text{ is continuous, with } c(0) = x \} \]
which we shall endow with the natural Wiener measure denoted $w_{t,x}$, and on which the Stratonovich integral $\operatorname{Strat} _{t,x}$ will live.

Let $V : M \to \mathbb R$ be continuous, with $\inf V > -\infty$ (in a future work it will be shown that the Feynman-Kac formula holds under much weaker conditions on $V$, at the price of significantly more complicated proofs; since the subject of this article is not the most general Feynman-Kac formula, but the functional-analytic and differential-geometric aspects of stochastic integration, the conditions imposed on $V$ above strike a satisfactory balance between clarity and generality). Consider the operator $(\mathrm d + \mathrm i \alpha)^* (\mathrm d + \mathrm i \alpha) + V$ acting on $C_0 ^\infty (M) \subset L^2 (M)$ (the star denotes the formal adjoint, i.e. the Hodge "$*$" operator); it is clearly symmetric and lower-bounded. Let $H_{\mathrm d + \mathrm i \alpha, V}$ be its densely-defined maximal self-adjoint extension obtained using the Friedrichs procedure.

\begin{theorem}[The Feynman-Kac-Itô formula]
If $f \in L^2(M)$, then
\[ (\mathrm e ^{-t H_{\mathrm d + \mathrm i \alpha, V}} f) (x) = \int _{\mathcal C _{t,x}} \mathrm e ^{\mathrm i \operatorname{Strat} _{t,x} (\alpha) - \int _0 ^t V(c(s)) \, \mathrm d s} \, f(c(t)) \, \mathrm d w_{t,x} (c) \]
for every $t>0$ and almost all $x \in M$.
\end{theorem}
\begin{proof}
Let us consider an exhaustion $M = \bigcup _{j \ge 0} U_j$ with regular domains, as already done previously in this work, the notations being the ones already encountered (essentially, every object related to $U_j$ gets an upper index ${}^{(j)}$). Using theorem 4 from \cite{Simon78}, we have that
\[ \mathrm e ^{-t H_{\mathrm d + \mathrm i \alpha, V}} = \lim _{j \to \infty} \mathrm e ^{-t H_{\mathrm d + \mathrm i \alpha, V} ^{(j)}} \]
strongly in $L^2(M)$, whereas from the Trotter-Kato formula we get that
\[ \mathrm e ^{-t H_{\mathrm d + \mathrm i \alpha, V} ^{(j)}} = \lim _{k \to \infty} \left( \mathrm e ^{-\frac t k H_{\mathrm d + \mathrm i \alpha, 0} ^{(j)}} \, \mathrm e ^{-\frac t k V} \right) ^k \]
strongly in $L^2(U_j)$, so that if $f, f' \in L^2(M)$, then
\begin{align*}
\langle (\mathrm e ^{-t H_{\mathrm d + \mathrm i \alpha, V}} & f), f' \rangle _{L^2(M)} = \lim _{j \to \infty} \langle (\mathrm e ^{-t H_{\mathrm d + \mathrm i \alpha, V} ^{(j)}} f), f' \rangle _{L^2 (U_j)} = \\
& = \lim _{j \to \infty} \lim _{k \to \infty} \left< \left[ \left( \mathrm e ^{-\frac t {2^k} H_{\mathrm d + \mathrm i \alpha, 0} ^{(j)}} \, \mathrm e ^{-\frac t {2^k} V} \right) ^{2^k} f \right], f' \right> _{L^2 (U_j)} = \\
& = \lim _{j \to \infty} \lim _{k \to \infty} \int _{U_j} \mathrm d x \left[ \int _{U_j} \mathrm d x_1 \, h_\alpha ^{(j)} \left( \frac t {2^k}, x, x_1 \right) \mathrm e ^{-\frac t {2^k} V(x_1)} \right. \dots \\
& \dots \left. \int _{U_j} \mathrm d x_{2^k} \, h_\alpha ^{(j)} \left( \frac t {2^k}, x_{2^k - 1}, x_{2^k} \right) \mathrm e ^{-\frac t {2^k} V(x_{2^k})} f(x_{2^k}) \right] \overline{f'(x)} = \\
& = \lim _{j \to \infty} \lim _{k \to \infty} \int _{U_j} \mathrm d x \left[ W_{\alpha, t, x} ^{(j)} \left( c \mapsto \mathrm e ^{- \frac t {2^k} \sum _{l=1} ^{2^k} V \left( c \left( \frac {tl} {2^k} \right) \right)} f(c(t)) \right) \right] \overline{f'(x)} = \\
& = \lim _{j \to \infty} \lim _{k \to \infty} \int _{U_j} \mathrm d x \left[ \int _{\mathcal C_{t,x} (\overline {U_j})} \rho_{\alpha, t, x} ^{(j)} (c) \, \mathrm e ^{- \frac t {2^k} \sum _{l=1} ^{2^k} V \left( c \left( \frac {tl} {2^k} \right) \right)} f(c(t)) \, \mathrm d w_{t,x} ^{(j)} (c) \right] \overline{f'(x)} = \\
& = \lim _{j \to \infty} \int _{U_j} \mathrm d x \left[ \int _{\mathcal C_{t,x} (\overline {U_j})} \rho_{\alpha, t, x} ^{(j)} (c) \, \mathrm e ^{- \int _0 ^t V(c(s)) \, \mathrm d s} f(c(t)) \, \mathrm d w_{t,x} ^{(j)} (c) \right] \overline{f'(x)} = \\
& = \int _M \mathrm d x \left[ \int _{\mathcal C_{t,x}} \rho_{\alpha, t, x} (c) \, \mathrm e ^{- \int _0 ^t V(c(s)) \, \mathrm d s} f(c(t)) \, \mathrm d w_{t,x} (c) \right] \overline{f'(x)} = \\
& = \int _M \mathrm d x \left[ \int _{\mathcal C_{t,x}} \mathrm e^{\mathrm i \operatorname{Strat} _{t,x} (\alpha)} \, \mathrm e ^{- \int _0 ^t V(c(s)) \, \mathrm d s} f(c(t)) \, \mathrm d w_{t,x} (c) \right] \overline{f'(x)} \ ,
\end{align*}
whence
\[ (\mathrm e ^{-t H_{\mathrm d + \mathrm i \alpha, V}} f) (x) = \int _{\mathcal C_{t,x}} \mathrm e^{\mathrm i \operatorname{Strat} _{t,x} (\alpha)} \, \mathrm e ^{- \int _0 ^t V(c(s)) \, \mathrm d s} f(c(t)) \, \mathrm d w_{t,x} (c) \]
for every $f \in L^2(M)$ and almost every $x \in M$. We have used the dominated convergence theorem twice; to see that its hypotheses are indeed met, notice that:
\begin{itemize}[wide]
\item $\rho_{\alpha, t, x} ^{(j)} \to \rho_{\alpha, t, x} = \mathrm e^{\mathrm i \operatorname{Strat} _{t,x} (\alpha)}$ almost everywhere on $\mathcal C_{t,x}$ by the very definitions of $\rho_{\alpha, t, x}$ and of $\operatorname{Strat} _{t,x} (\alpha)$;
\item the functions $\rho_{\alpha, t, x} ^{(j)}$ and $\rho_{\alpha, t, x}$ are bounded in absolute value by $1$ almost everywhere;
\item $\frac t {2^k} \sum _{l=1} ^{2^k} V \left( c \left( \frac {tl} {2^k} \right) \right) \to \int _0 ^t V(c(s)) \, \mathrm d s$, for every continuous curve $c$ (because $V$ has been assumed continuous), so the function $c \mapsto \int _0 ^t V(c(s)) \, \mathrm d s$ is measurable as the pointwise limit of a sequence of measurable functions;
\item the functions $\mathrm e ^{- \frac t {2^k} \sum _{l=1} ^{2^k} V \left( c \left( \frac {tl} {2^k} \right) \right)}$ and $\mathrm e ^{- \int _0 ^t V(c(s)) \, \mathrm d s}$ are bounded by $\mathrm e ^{-t \inf V}$;
\item the function $c \mapsto f(c(t))$ is integrable on $\mathcal C _{t,x}$ because
\[ \int _{\mathcal C _{t,x}} |f(c(t))| \, \mathrm d w_{t,x} (c) = \int _M h(t,x,y) \, |f(y)| \, \mathrm d y \le \sqrt {\int _M h(t,x,y)^2 \, \mathrm d y \int _M |f(y)| ^2 \, \mathrm d y} \ ; \]
(a similar argument is used to show the integrability on $\mathcal C _{t,x} (\overline {U_j})$);
\item the Wiener measures $w_{t,x} ^{(j)}$ and $w_{t,x}$ are finite.
\end{itemize}
\end{proof}

\section{Conclusion and acknowledgements}

The main aim has been to show how to give an alternative construction of some basic objects in stochastic analysis using only functional-analytic tools, without it being necessary to resort to probability-theoretical concepts or techniques. Another aim has been to advance a point of view allowing the entire subject of stochastic integration to be seen unravelling from a small number of fundamental ideas, along lines emphasizing the deep analogies with curvilinear integration. The strategy adopted herein has allowed the classification of stochastic integrals and the displaying of the simple relationship connecting any two of them. In particular, we have seen that the Stratonovich integral is the natural stochastic instrument in problems with a strong geometric flavour, it allowing for many differential-geometric ideas to be transported almost unchanged into the stochastic framework. A future work, currently in draft form, will show that the Itô integral is the appropriate tool in those problems with a strong probabilistic (or analytic) flavour. None of these two integrals is "better" than the other one, the choice between them being made in connection with the specifics of the problem under study.

From a technical point of view, since all the objects involved were intrinsic to the manifold $M$, we have obtained that their construction be intrinsic, too. This differs from the approach that other stochastic analysis on manifolds textbooks use (for instance \cite{Hsu02}), which resort to embedding the underlying manifold in Euclidean spaces using Whitney's theorem, thus using extrinsic geometrical tools to obtain intrinsic results. We have also attempted to keep the prerequisites to a minimum, using only a handful of basic functional-analytic tools. Once the foundations of this construction are laid down, developing the various properties of stochastic integrals becomes much easier than in the traditional probability-theoretic textbooks. We have thus not needed to use Cartan's rolling map, as it is done in \cite{Elworthy82}. Neither has it been necessary to choose an interpolation rule, as done in \cite{Emery89} (which requires the use of the measurable selection theorem, checking the hypotheses of which further requires working with the Whitney topology on the space of smooth curves in $M$), its role being taken on by the cut-off function $\chi$ as well as by the truncation by $0$ of the expression $\int _{[0,1]} \alpha _{\gamma_{x,y} (\tau)} (\dot \gamma_{x,y} (\tau)) \, \mathrm d P (\tau)$ for $y$ far away from $x$. Unlike in \cite{Duncan76}, $M$ is not required to be compact. We have also not needed to work with second order tangent vectors and Laurent Schwartz' second order differential geometry, as done by Émery in \cite{Emery89}. This parsimonious use of fundamental concepts and technical means has been one of the driving goals of the present text which is built upon the belief that conceptual and technical minimality must be an imperative of any intellectual construction.

The key points to remember from this article are:
\begin{itemize}[wide]
\item for every regular Borel probability $P$ on the interval $[0,1]$ there exists a unique corresponding stochastic integral $\operatorname{Int}_{P,t} : \Omega ^1 (M) \to L^0 (\mathcal C_t)$, that is linear;
\item if $P$ and $Q$ have the same first-order moment, then $\operatorname{Int}_{P,t} = \operatorname{Int}_{Q,t}$;
\item $\operatorname{Int}_{\delta_0,t}$ is the Itô integral;
\item if $P = \operatorname{Leb}_{[0,1]}$, or $P = \delta_{\frac 1 2}$, or $P = \frac 1 2 (\delta_0 + \delta_1)$, then the corresponding stochastic integral $\operatorname{Int}_{P,t}$ is the Stratonovich integral;
\item for any probabilities $P$ and $Q$ as above,
\[ \operatorname{Int}_{P,t} (\alpha) (c) = \operatorname{Int}_{Q,t} (\alpha) (c) - 2 \int _{[0,1]} \tau \, \mathrm d (P-Q) (\tau) \int _0 ^t (\mathrm d^* \alpha) (c(s)) \, \mathrm d s \ ; \]
\item $\operatorname{Int}_{P,t} (\alpha)$ is the limit in measure of the approximations $c \mapsto \sum _{j=0} ^{2^k - 1} I_P (\alpha) \left( c \left( \frac {jt} {2^k} \right), c \left( \frac {(j+1)t} {2^k} \right) \right)$.
\end{itemize}

That the stochastic integrals, as defined in this work, form a sort of continuum, having the Itô and Stratonovich ones among them as particular cases, is an idea that has been explored in \cite{PLS13}, too, but only for $M = \mathbb R$, from a completely different point of view and using entirely different mathematical tools.

Reaching the end of this work, it is a pleasure to thank Dr. Radu Purice from the "Simion Stoilow" Institute of Mathematics of the Romanian Academy for his unabated support, both mathematical and moral, generously offered to me during the difficult elaboration of the present work. His company during the hesitant explorations of these uncharted waters, his constant availability for endless discussions in the margin of this text and his patient reading of its many preliminary versions are the invisible and discrete ingredients without which this article would have remained just a scribble.

\end{document}